\documentclass{article}
\usepackage{graphicx} 
\usepackage{amsfonts}
\usepackage{bbm}
\usepackage{amscd}
\usepackage{xypic}
\usepackage{amsmath}
\usepackage{amssymb}
\usepackage{amsthm}
\usepackage{bbm}
\usepackage{CJK}
\usepackage{fancyhdr}
\usepackage{hyperref}
\usepackage{indentfirst}
\usepackage{latexsym}
\usepackage{mathrsfs}
\allowdisplaybreaks 
\usepackage{color}
\usepackage{amsmath}
\numberwithin{equation}{section}
\usepackage[top=1in,bottom=1in,left=1.25in,right=1.25in]{geometry}
\title{A classification of pre-Lie $H$-pseudoalgebras of low ranks}
\author{Botong Gai}
\date{October 2025}

\begin{document}

\author {{\bf Botong Gai}$^a$\\
{\small a: School of Mathematics, Southeast University,  Nanjing, Jiangsu}\\
{\small 210096, P. R. of China. E-mail: 230228425@seu.edu.cn}\\
{\small Nanjing, Jiangsu 210096, P. R. of China}}
\maketitle
\begin{center}
		\begin{minipage}{14.cm}
               {\bf Abstract} Let $H=U(\delta)$ be the universal enveloping algebra of finite dimension Lie algebra $\delta$.  The central result of the paper is the classification of pre-Lie $H$-pseudoalgebras of low ranks over the Hopf algebra $H$. We firstly study pre-Lie pseudoalgebras that are free of rank $1$ over $H$. Then we introduce and classify a class of pre-Lie $H$-pseudoalgebras $\mathcal{P}$ which are generated by two pre-Lie pseudoalgebras of rank $1$. Finally, the associativity of $\mathcal{P}$ is also considered and a explicit assification is presented.
\\

\\


		\end{minipage}
	\end{center}

\section{Introduction}
\def\theequation{0. \arabic{equation}}
	\setcounter{equation} {0} \hskip\parindent

The notion of conformal algebras (\cite{K}) was introduced by Kac as an axiomatic description of
the operator product expansion (OPE) of chiral fields in conformal field theory, and it came to
be useful for investigation of vertex algebras. In 2001, Bakalov, D'Andrea and Kac introduced a notion of “multi-dimensional” Lie conformal algebras by replacing the polynomial algebra $\mathbb{C}[\partial]$
  with any cocommutative Hopf algebra $H$, called Lie $H$-pseudoalgebras (see \cite{BDK}). Furthermore, they established such theory as
  cohomology theory, representation theory of $H$-pseudoalgebras
   and  irreducible modules over finite simple Lie pseudoalgebras (see \cite{BDK}-\cite{BDK3}, \cite{DL}). Actually,
  Lie $H$-pseudoalgebras  can be considered as Lie algebras in a certain ``pseudotensor'' category, instead of the category
of vector spaces. A pseudotensor
category (\cite{BD}) is a category equipped with ``polylinear maps'' and a way to compose them
(such categories were first introduced by Lambek (\cite{L}) under the name multi-categories). This is enough to define the notions of Lie algebra, representations, cohomology, etc.
 \\

The classification of pseudoalgebras is one of important subjects that scholars are interested in. It is well-known that the  Lie conformal algebra of rank $1$ is either abelian Lie conformal
algebra or Virasoro Lie conformal algebra (\cite{DK}). The classification Leibniz conformal algebras
 of rank $1$ and rank $2$  was considered in  pseudoalgebra method (\cite{W}). Besides, the Lie conformal algebras of rank $2$
are also studied by many other methods (see, \cite{BCH}, \cite{HL}). In \cite{BDK}, the authors classified the simple Lie $U(\delta)$-pseudoalgebra of rank $1$   into two types: $H(\delta, \chi, \omega)$ and $K(\delta, \theta)$, which is defined as Lie $H$-pseudoalgebra analogues of the primitive linearly compact Lie algebras $H_{N}$ and $K_{N}$. Also in this paper, the classification of simple and semisimple Lie $U(\delta)$-pseudoalgebras are finished. Based on these results, a class Lie pseudoalgebras of rank $2$  which is generated by an abelian Lie $H$-pseudoalgebra of rank $1$ and a simple Lie $H$-pseudoalgebra was characterized and studied in \cite{GW}.
\\

The notion of pre-Lie algebra (also called left-symmetric algebras, quasi-associative algebras, Vinberg algebras and
so on) is introduced independently by M. Gerstenhaber in deformation theory of rings and algebras (\cite{G}). Pre-Lie
algebra arose from the study of affine manifolds and affine structures on Lie group (\cite{KO}), homogeneous convex cones (\cite{V}).
Its defining identity is weaker than associativity. This algebraic structure describes some properties of cochains space in
Hochschild cohomology of an associative algebra, rooted trees and vector fields on affine spaces. Moreover, it is playing
an increasing role in algebra, geometry and physics due to their applications in nonassociative algebras, combinatorics,
numerical analysis and quantum field theory (see \cite{BA1}-\cite{BA2},\cite{BU}). In \cite{LW}, the authors introduce pre-Lie algebra in pseudotensor category, which is called pre-Lie pseudoalgebra. They obtain a large number of pre-Lie $H$-pseudoalgebras by using Rota-Baxter operators and study their annihilation algebras.\\

The organization of this paper is as follows. In Section 1, we recall some necessary definitions, notations and some results about Hopf algebras and pre-Lie $H$-pseudoalgebras. In Section 2, we study the (left) pre-Lie $H$-pseudoalgebra of rank one and give a basic characterization. In Section 3, a class of pre-Lie $H$-pseudoalgebras $\mathcal{P}$ are introduced, which are generated by two rank one pre-Lie $H$-pseudoalgebras. Further, we present complete classifications of $\mathcal{P}$, which can actually be seen as a method to constructing pre-Lie pseudoalgebra. In Section 4, we discuss the associativity of $\mathcal{P}$ and present a simple classification.\\

\section*{1. Preliminaries}
	\def\theequation{1. \arabic{equation}}
	\setcounter{equation} {0} \hskip\parindent

 Unless otherwise specified, all vector spaces, linear maps and tensor products are considered over an algebraically
 closed field $\mathbf{k}$ of characteristic $0$.

 \subsection*{1.1 Hopf algebra}

In this section we present some facts and notations of Hopf algebras which will be used throughout this paper and can be found, for example, in Sweedler's book (cf. \cite{S}). Let $H$ be a Hopf algebra with a coproduct $\Delta$ , a counit $\epsilon$, and an antipode $S$.
We will use the standard Sweedler's notation (cf. \cite{S}): $\Delta(h) = h_{(1)}\otimes h_{(2)}$. The axioms of the antipode and the counit
can be written as follows:

\begin{equation*}
   S(h_{(1)})h_{(2)}=\epsilon(h)1_{H}=h_{(1)}S(h_{(2)}), \quad \epsilon(h_{(1)})h_{(2)}=h_{(1)}\epsilon(h_{(2)})=h,
\end{equation*}

while the fact that $\Delta$ is an algebra homomorphism translates as
  $(hg)_{(1)} \otimes(hg)_{(2)}=h_{(1)}g_{(1)} \otimes h_{(2)}g_{(2)}$. An element $x\in H$ is said to be a group-like element if $\Delta (x)=x\otimes x$.
   We denote $G(H)$ by the set of all the group-like elements of $H$.
\\

We consider an increasing sequence of subspaces of a Hopf algebra $H$ defined inductively by:
\begin{align*}
  &F^{0}H=\mathbf{k}[G(H)],  \quad \quad F^{n}H=0, \quad for \quad n<0,
\end{align*}
\begin{center}
  $F^{n}H=span_{\mathbf{k}}\{h\in H|\Delta(h)\in F^{0}H\otimes h+h\otimes F^{0}H+\sum_{i=1}^{n-1}F^{i}H\otimes F^{n-i}H\}$, for $n\geq 1$.
\end{center}
It has the following properties (which are immediate from definitions):
\begin{equation*}
  (F^{n}H)(F^{m}H)\subset F^{n+m}H, \quad \Delta(F^{n}H)\subset\sum_{i=0}^{n}F^{i}H\otimes F^{n-i}H, \quad S(F^{n}H)\subset F^{n}H.
\end{equation*}
When $H$ is cocommutative, one can show that:
\begin{equation}
  \bigcup_{n} F^{n}H=H.
\end{equation}
(This condition is also satisfied when $H$ is a quantum universal enveloping algebra.) At this time, we say that a nonzero element $a\in H$ has \textbf{degree $n$} if $a \in F^{n} H \setminus F^{n-1}H$.
\\

When $H$ is an universal enveloping algebra of a finite-dimensional Lie algebra, or its
smash product with the group algebra of a finite group, the following finiteness condition holds:
\begin{equation}
  dim F^{n}H<\infty.
\end{equation}

Let $\delta$ be a finite-dimensional Lie algebra over $\mathbf{k}$. Let $\{\partial_{i},...,\partial_{N}\}$ be a basis of $\delta$ and $\partial^{I}=\frac{\partial_{1}^{i_{1}}...\partial_{N}^{i_{N}}}{i_{1}!...i_{N}!}$ for $I=(i_{1},i_{2},...,i_{N})\in \mathbb{Z}_{+}^{N}$. Then $\{\partial^{I}\}$ is a basis of $H$ (the \emph{Poincare-Birkhoff-Witt (PBW)} basis). Moreover, the coproduct on $\partial^{I}$ is
\begin{equation}
  \Delta(\partial^{I})=\sum_{J+K=I}\partial^{J}\otimes\partial^{K}.
\end{equation}
In the following, we set $|I|=i_{1}+i_{2}+...+i_{N}$ for convenience. Further, if $a=\partial_{1}^{j_{1}}...\partial_{N}^{j_{N}}\in U(\delta)$, then the degree of $a$ denotes by $deg(a)$ is actually $j_{1}+j_{2}+...+j_{N}$.
\\

 For an arbitrary Hopf algebra $H$,  a map $\mathcal{F} : H \otimes H \rightarrow H \otimes H$ is called the
\emph{Fourier transform} (or \emph{Galois map}), by the formula, for all $f, g\in H$

\begin{equation*}
  \mathcal{F}(f\otimes g)=fS(g_{(1)})\otimes g_{(2)}.
\end{equation*}

Observe that $\mathcal{F}$ is a vector space isomorphism with an inverse given by
\begin{center}
  $\mathcal{F}^{-1}(f\otimes g)=fg_{(1)}\otimes g_{(2)}$.
\end{center}
The significance of $\mathcal{F}$ is in the identity

\begin{equation*}
  f\otimes g=\mathcal{F}^{-1}\mathcal{F}(f\otimes g)=(fS(g_{(1)})\otimes 1)\Delta(g_{(2)}),
\end{equation*}

which implies the following result.
\\

\textbf{Lemma 1.1.}(\cite{BDK} or \cite{VD1}) Every element of $H \otimes H$ can be uniquely represented in the form $\sum_{i}(h_{i}\otimes 1)\Delta(l_{i})$, where ${h_{i}}$ is a fixed $\mathbf{k}$-basis of $H$ and $l_{i}\in H$. In other words, $H\otimes H = (H \otimes \mathbf{k})\Delta(H)$.
\\

\subsection*{1.2. $H$-pseudoalgebra}

In this section, we shall recall the basic definitions and examples of $H$-pseudoalgebra.
\\

 An $H$-pseudoalgebra $(A, \ast)$ is a left $H$-module $A$ together with a map
(called the pseudo-product):
\begin{equation*}
  \ast: A\otimes A\rightarrow (H\otimes H)\otimes_{H} A, \qquad x\otimes y\mapsto x\ast y,
\end{equation*}
satisfying \textbf{$H$-bilinearity}: for any $x, y\in A, h, g \in H$,
\begin{equation}
hx\ast gy=(h\otimes g\otimes_{H} 1)(x\ast y).
\end{equation}
Specifically, if $x\ast y=\sum_{i}h_{i}\otimes g_{i}\otimes_{H}e_{i}$, then $hx\ast gy=\sum_{i}hh_{i}\otimes gg_{i}\otimes_{H}e_{i}$.

Further, we call the $H$-pseudoalgebra $(A, \ast)$ is \textbf{associative} if
\begin{equation}
  (x\ast y)\ast z=x\ast(y\ast z)
\end{equation}
for any $x, y, z\in A$.

We call the $H$-pseudoalgebra $(A, \ast)$ is \textbf{left pre-Lie} (respectively, \textbf{right pre-Lie}) if
\begin{align}
  (x\ast y)\ast z-x\ast(y\ast z)=&((12)\otimes_{H} \mathrm{id})[(y\ast x)\ast z-y\ast(x\ast z)] \\
  \bigl(resp.(x\ast y)\ast z-x\ast(y\ast z)=&((23)\otimes_{H} \mathrm{id})[(x\ast z)\ast y-x\ast(z\ast y)] \bigr)
\end{align}
in $H^{\otimes3}\otimes_{H} A$, where $(12)(f\otimes h\otimes g)=h\otimes f\otimes g$, $(23)(f\otimes h\otimes g)=f\otimes g\otimes h$ and for $\beta\in H\otimes H$,
\begin{align*}
  &(\beta\otimes_{H} x)\ast y= \sum_{i} (\beta\otimes 1)(\Delta\otimes \mathrm{id})(\alpha_{i})\otimes_{H}c_{i},\\
  &x\ast (\beta\otimes_{H} y)= \sum_{i} (1\otimes \beta)(\mathrm{id}\otimes\Delta)(\alpha_{i})\otimes_{H}c_{i},
\end{align*}
if $x\ast y=\sum_{i} \alpha_{i}\otimes_{H} c_{i}\in H\otimes H\otimes_{H} A$.\\

\textbf{Remark.} Note that an associative $H$-pseudoalgebra is naturally a pre-Lie $H$-pseudoalgebra while the opposite is usually not true.
\\

Combine with Lemma 1.1, we need to point out that the pseudo-product of pre-Lie $H$-pseudoalgebras is well defined, provided that the Hopf algebra $H$ is cocommutative. We will always assume that $H$ is cocommutative when talking about pre-Lie $H$-pseudoalgebras.\\

\textbf{Example.} Let $H'$ be a Hopf subalgebra of $H$ and $(A, \ast)$ be a pre-Lie $H'$-pseudoalgebra. Then we define the current $H$-pseudoalgebra $Cur^{H}_{H'} A \equiv Cur A$ as $H\otimes_{H'} A$ by extending the
pseudo-product $a \ast b$ of $A$ by $H$-bilinearity.  Explicitly, for $a, b\in A$, $f, g\in H$,
\begin{center}
  $(f\otimes_{H'} a) \ast (g\otimes_{H'} b)=(f\otimes g\otimes_{H} 1)(a\ast b)$.
\end{center}

A special case is when $H'=\mathbf{k}$: given a pre-Lie algebra $(A, \circ)$, let
$Cur A=H\otimes A$ with the following pseudo-product:
\begin{center}
  $(f\otimes_{H'} a) \ast (g\otimes_{H'} b)=(f\otimes g)\otimes_{H} (1\otimes a\circ b)$.
\end{center}
Then $Cur A$ is a pre-Lie $H$-pseudoalgebra.
\\

\section*{2.  Pre-Lie pseudoalgebras of rank one}
	\def\theequation{2. \arabic{equation}}
	\setcounter{equation} {0} \hskip\parindent

Unless otherwise specified, we will often be working with the Hopf algebra $H=U(\delta)$ in the following, where $\delta$ is a
finite-dimensional Lie algebra over $\mathbf{k}$.\\

Let $A=He$ be a left pre-Lie $H$-pseudoalgebra which is a free $H$-module of rank 1. Then by $H$-bilinearity, the pseudo-product on $A$ is determined by $e\ast e$, or equivalently, by an $\alpha\in H\otimes H$ such that $e\ast e=\alpha\otimes_{H} e$.\\

\textbf{Proposition 2.1.} $A=He$ with the pseudo-product $e\ast e=\alpha\otimes_{H} e$ is a left (resp. right) pre-Lie $H$-pseudoalgebra if and only if $\alpha\in H\otimes H$ satisfying the following equation:
\begin{align}
  &(\alpha\otimes1)(\Delta\otimes \mathrm{id})\alpha-(1\otimes\alpha)(\mathrm{id}\otimes\Delta)\alpha=(\sigma\otimes \mathrm{id})[(\alpha\otimes1)(\Delta\otimes \mathrm{id})\alpha-(1\otimes\alpha)(\mathrm{id}\otimes\Delta)\alpha].\\
\bigl(resp. &(\alpha\otimes1)(\Delta\otimes \mathrm{id})\alpha-(1\otimes\alpha)(\mathrm{id}\otimes\Delta)\alpha=(\mathrm{id}\otimes \sigma)[(\alpha\otimes1)(\Delta\otimes \mathrm{id})\alpha-(1\otimes\alpha)(\mathrm{id}\otimes\Delta)\alpha]. \bigr)
\end{align}

\textbf{Proof.} This follows immediately from (1.6) and (1.7).
$\hfill \blacksquare$
\\

\textbf{Proposition 2.2.} Let $H$ be an universal enveloping algebra of Lie
algebra $\delta$. Then the solution $\alpha\in H\otimes H$ of equation (2.1) is $\alpha=1\otimes s+t\otimes1$ for arbitrary $s\in\delta, t\in\mathbf{k}$.

\textbf{Proof.}
Obviously, $\alpha=0$ is a solution of (2.1). If $\alpha\neq0$ and $\alpha=\sigma(\alpha)$, then we can set $\alpha=\sum_{i}a_{i}\otimes b_{i}+b_{i}\otimes a_{i}$ for some nonzero $\{a_{i}, b_{i}\}\subseteq H$. One can easily verify
\begin{equation*}
  (\alpha\otimes1)(\Delta\otimes \mathrm{id})\alpha=(\sigma\otimes\mathrm{id})[(\alpha\otimes1)(\Delta\otimes \mathrm{id})\alpha]
\end{equation*}
 and (2.1) is actually
\begin{equation}
  (1\otimes\alpha)(\mathrm{id}\otimes\Delta)\alpha=(\sigma\otimes\mathrm{id})[(1\otimes\alpha)(\mathrm{id}\otimes\Delta)\alpha].
\end{equation}
We have proved in \cite{GW} (Lemma 2.4) that solutions of (2.3) are the form of $h\otimes1$ for $h\in H$. Since $\alpha=\sigma(\alpha)$, there must have $\alpha=t\otimes1$ for nonzero $t\in\mathbf{k}$.

If $\alpha\neq0$ and $\alpha\neq\sigma(\alpha)$. Let $\{\partial_{1}, ..., \partial_{N}\}$ be a basis of $\delta$. We consider the corresponding
\emph{PBW} basis of $H=U(\delta)$ given by $\partial^{I}=\partial_{1}^{i_{1}}...\partial_{N}^{i_{n}}/i_{1}!...i_{N}!$, where $I=(i_{1},...,i_{N})\in \mathbb{Z}_{+}^{N}$. Then we can write $\alpha=\sum_{I}\alpha_{I}\otimes\partial^{I}$, $\alpha_{I}\in H$ and equation (2.1) becomes
\begin{equation*}
  \sum_{I}[\alpha\Delta(\alpha_{I})-\sigma\bigl(\alpha\Delta(\alpha_{I})\bigr)]\otimes\partial^{I}=\sum_{I,J,K}(\alpha_{J+K}\otimes\alpha_{I}\partial^{J}-\alpha_{I}\partial^{J}\otimes\alpha_{J+K})\otimes\partial^{I}\partial^{K}.
\end{equation*}
Notice that the left summation is not equal to zero, then similar to Lemma 4.1 in \cite{BDK}, we have $|I|\leq1$ and $\alpha=\sum_{i=1}^{N}h_{i}\otimes\partial_{i}+g\otimes1+1\otimes s+t\otimes1$ for some $h_{i}, g\in H, s\in\delta, t\in\mathbf{k}$ where $h_{i}, g$ do not contain constant terms. Taking this into (2.1) and comparing the degree of third tensor factor, we obtain
\begin{multline}
    \sum_{i,j}(h_{i}\otimes\partial_{i}+g\otimes1+1\otimes s-\partial_{i}\otimes h_{i}-1\otimes g-s\otimes1)\Delta(h_{j})\otimes\partial_{j}+(g+s)\otimes1\otimes s
-1\otimes (g+s)\otimes s\\
    =\sum_{i,j}(h_{j}\otimes h_{i}\partial_{j}-h_{i}\partial_{j}\otimes h_{j})\otimes\partial_{i}+\sum_{i,j}h_{j}\otimes h_{i}\otimes[\partial_{i}, \partial_{j}]
+\sum_{i}(1\otimes h_{i}-h_{i}\otimes1)\otimes[\partial_{i}, s]\\
    +\sum_{i}(1\otimes h_{i}s-h_{i}s\otimes1)\otimes\partial_{i},
\end{multline}
and
\begin{multline}
  \sum_{i}(h_{i}\otimes\partial_{i}+g\otimes1+1\otimes s-\partial_{i}\otimes h_{i}-1\otimes g-s\otimes1)\Delta(g)+t(g\otimes1-1\otimes g)\\=\sum_{i}h_{i}\otimes g\partial_{i}-g\partial_{i}\otimes h_{i}+1\otimes gs-gs\otimes1.
\end{multline}
Applying $\mathrm{id}\otimes\epsilon$ to both sides of (2.5), we obtain $(g+t)g=[s,g]$.  Since $g$ does not contain constant term, if $g\neq0$, then $deg(g)=p\geq1$ and the degree of left hand is exactly $2p$ while the right is $p$ at most, which gives a contradiction. Thus, there must have $g=0$
and equation (2.4) becomes
\begin{multline}
    \sum_{i,j}(h_{j}\otimes h_{i}\partial_{j}-h_{i}\partial_{j}\otimes h_{j})\otimes\partial_{i}+\sum_{i,j}h_{j}\otimes h_{i}\otimes[\partial_{i}, \partial_{j}]+\sum_{i}(1\otimes h_{i}-h_{i}\otimes1)\otimes[\partial_{i}, s]
    \\+\sum_{i}(1\otimes h_{i}s-h_{i}s\otimes1)\otimes\partial_{i}= \sum_{i,j}(h_{i}\otimes\partial_{i}+1\otimes s-\partial_{i}\otimes h_{i}-s\otimes1)\Delta(h_{j})\otimes\partial_{j}.
\end{multline}
Suppose there exists $h_{k}\neq0$ among all $h_{i}$'s, then there will exist
some (nonzero) $h_{j}$ of maximal degree $d$. Since $h_{i}$'s do not contain constant term, we have $d\geq1$. If $d>1$, then $2d>d+1$. There exists term in the right hand lies in $F^{2d}(H)\otimes F^{1}(H)\otimes F^{1}(H)$, which can not be cancelled by other terms. Thus, $d=1$, $h_{i}\in\delta$ and equation (2.6) becomes
\begin{multline*}
  \sum_{i,j}[h_{i}h_{j}\otimes\partial_{i}-[\partial_{i}, h_{j}]\otimes h_{i}+h_{j}\otimes s-sh_{j}\otimes1+h_{i}\otimes[\partial_{i}, h_{j}]-\partial_{i}\otimes h_{i}h_{j}+1\otimes sh_{j}-s\otimes h_{j}]\otimes\partial_{j}\\
  =\sum_{i,j}h_{j}\otimes h_{i}\otimes[\partial_{i}, \partial_{j}]+\sum_{i}(1\otimes h_{i}-h_{i}\otimes1)\otimes[\partial_{i}, s]+\sum_{i}(1\otimes h_{i}s-h_{i}s\otimes1)\otimes\partial_{i}.
\end{multline*}
Notice that there exist terms $\sum_{i,j}h_{i}h_{j}\otimes\partial_{i}\otimes\partial_{j}$ lie in $F^{2}(H)\otimes F^{1}(H)\otimes F^{1}(H)$, which can not be cancelled by other terms and then gives a contradiction. Therefore, we have $h_{i}=0~(\forall i)$ and $\alpha=1\otimes s+t\otimes1$ for some nonzero $s\in\delta$.

Summarizing discussion above, we have $\alpha=1\otimes s+t\otimes1$ for arbitrary $s\in\delta, t\in\mathbf{k}$.
$\hfill \blacksquare$
\\

The proposition above points out that $He$ is rank 1 left pre-Lie $H$-pseudoalgebra if and only if $e\ast e=(1\otimes s+t\otimes1)\otimes_{H}e$, where $s\in\delta, t\in\mathbf{k}$ are arbitrary.
Similarly, one can also prove the following result:\\

\textbf{Proposition 2.3.} Let $H$ be an universal enveloping algebra of Lie
algebra $\delta$. Then the solution $\alpha\in H\otimes H$ of equation (2.2) is $\alpha=s\otimes 1+t\otimes1$ for arbitrary $s\in\delta, t\in\mathbf{k}$, which means $He$ is rank 1 right pre-Lie $H$-pseudoalgebra if and only if $e\ast e=(s\otimes 1+t\otimes1)\otimes_{H}e$, where $s\in\delta, t\in\mathbf{k}$ are arbitrary.
\\

\textbf{Corollary 2.4.} Let $A=He$ be a rank one pre-Lie $H$-pseudoalgebra (both left and right). Then $e\ast e=t\otimes1\otimes_{H} e$ for arbitrary $t\in\mathbf{k}$.
\\

\section*{3.  Pre-Lie pseudoalgebras of rank two}
	\def\theequation{3. \arabic{equation}}
	\setcounter{equation} {0} \hskip\parindent

In this section, we consider a semidirect product $\mathcal{P}=He_{1}\oplus He_{2}$ of rank two, which is generated by $e_{1}$, $e_{2}$ as $H$-module. We discuss the structures of left pre-Lie $H$-pseudoalgebra on $\mathcal{P}$ and present classifications in all cases, which also can be seen as a method of constructing pre-Lie $H$-pseudoalgebras.\\

We firstly present some Lemmas, which are useful in the following.\\

\textbf{Lemma 3.1.} (\cite{GW}) Solution $\alpha \in H\otimes H$ of the following equation
\begin{equation}
  (1\otimes\alpha)(\mathrm{id}\otimes \Delta)\alpha=(\sigma\otimes \mathrm{id})[(1\otimes\alpha)(\mathrm{id}\otimes \Delta)\alpha]
\end{equation}
is the form $\alpha=h\otimes 1$ for $h\in H$.
\\

\textbf{Lemma 3.2.} Let $s\in \delta$, $t\in \mathbf{k}$ and $s\neq0$. Then the solution $\alpha \in H\otimes H$ of the following equation
\begin{equation}
  (1\otimes s\otimes1+t)(\Delta\otimes \mathrm{id})\alpha-(1\otimes\alpha)(\mathrm{id}\otimes\Delta)\alpha=(\sigma\otimes \mathrm{id})[(1\otimes s\otimes1+t)(\Delta\otimes \mathrm{id})\alpha-(1\otimes\alpha)(\mathrm{id}\otimes\Delta)\alpha]
\end{equation}
is $\alpha=0$ or $\alpha=1\otimes s+ls\otimes1+k$ for arbitrary $l, k\in\mathbf{k}$.

\textbf{Proof.} Let $\{\partial_{1}, ..., \partial_{N}\}$ be a basis of $\delta$ and we also consider the corresponding
\emph{PBW} basis given by $\partial^{I}=\partial_{1}^{i_{1}}...\partial_{N}^{i_{n}}/i_{1}!...i_{N}!$.
We set $\alpha=\sum_{I}\alpha_{I}\otimes\partial^{I}$ and then (3.2) becomes
\begin{equation}
  \sum_{I}(1\otimes s-s\otimes1)\Delta(\alpha_{I})\otimes\partial^{I}=\sum_{I,J,K}(\alpha_{J+K}\otimes\alpha_{I}\partial^{J}-\alpha_{I}\partial^{J}\otimes\alpha_{J+K})\otimes\partial^{I}\partial^{K}.
\end{equation}
Similar to Lemma 4.1 in \cite{BDK}, we know $|I|\leq1$. Then we can write $\alpha=\sum_{i}\alpha_{i}\otimes\partial_{i}+g\otimes1$ for some $\alpha_{i}$'s, $g\in H$. Taking this into (3.3) and comparing the degree of the third tensor factor, we obtain
\begin{align}
  (1\otimes s-s\otimes1)\Delta(g) &= \sum_{i}\alpha_{i}\otimes g\partial_{i}-g\partial_{i}\otimes\alpha_{i}, \\
  \sum_{i}(1\otimes s-s\otimes1)\Delta(\alpha_{i})\otimes\partial_{i} &= \sum_{i,j}(\alpha_{j}\otimes\alpha_{i}\partial_{j}-\alpha_{i}\partial_{j}\otimes\alpha_{j})\otimes\partial_{i}+\alpha_{j}\otimes\alpha_{i}\otimes[\partial_{i},\partial_{j}].
\end{align}
Let $d$ be the maximal value among all $deg(\alpha_{i})$'s for $\alpha_{i}\neq0$, then terms
\begin{equation*}
\sum_{deg(\alpha_{k})=d}\alpha_{k}\partial_{k}\otimes\alpha_{k}\otimes\partial_{k}
\end{equation*}
in the right of (3.5) lie in $F^{d+1}(H)\otimes F^{d}(H)\otimes\delta$ while the only other terms that the first tensor factor lies in $F^{d+1}(H)$ are
\begin{equation*}
  \sum_{deg(\alpha_{k})=d}s\alpha_{k}\otimes1\otimes\partial_{k}.
\end{equation*}
Then there must have $d=0$ and we can rewrite $\alpha=1\otimes a+h\otimes1+k\otimes1$ for some $a\in\delta, k\in\mathbf{k}, h\in H$, where $h$ does not contain constant term. Taking this into (3.2), we obtain
\begin{multline}
  1\otimes s\otimes a+(1\otimes s)\Delta(h)\otimes1+1\otimes ks\otimes1-1\otimes a\otimes a-1\otimes ha\otimes1-1\otimes ka\otimes1\\
=s\otimes1\otimes a+(s\otimes1)\Delta(h)\otimes1+ks\otimes1\otimes1-a\otimes1\otimes a-ha\otimes1\otimes1-ka\otimes1\otimes1
\end{multline}
Applying $\mathrm{id}\otimes\epsilon\otimes\epsilon$ to both sides, we get $k(a-s)=sh-ha=[s, h]+h(s-a)$, which induces two cases as follows:

\textbf{Case I.} If $h=0$, one can easily obtain
 \begin{equation*}
  \begin{cases}
k = 0 \\
a\in\delta ~is~arbitrary
\end{cases}
\quad or \quad
\begin{cases}
k\in\mathbf{k}~is~arbitrary \\
a=s
\end{cases}.
 \end{equation*}
The second solution will make (3.6) holds all the time. Taking the first solution into (3.6), one can obtain
\begin{equation*}
1\otimes s\otimes a-1\otimes a\otimes a=s\otimes1\otimes a-a\otimes1\otimes a,
\end{equation*}
which implies $a=0$ or $a=s$. Therefore, $\alpha=0$ or $\alpha=1\otimes s$ or $\alpha=1\otimes s+k$ where $k$ is arbitrary.

\textbf{Case II.} If $h\neq0$, since $h\in H$ and $h$ does not contain constant term, we have $deg(h)=p\geq1$. Suppose $a\neq s$, then $deg[h(s-a)]>deg(h)\geq1$ while $deg[s, h]\leq deg(h)$ and $deg[k(a-s)]\leq1$, which is a contradiction. Thus, there must have $a=s$ and (3.6) becomes
\begin{equation*}
  (1\otimes s)\Delta(h)-1\otimes hs=(s\otimes1)\Delta(h)-hs\otimes1.
\end{equation*}
By comparing the coefficients of intermediate corresponding degree terms of $(1\otimes s)\Delta(h)$ and $(s\otimes1)\Delta(h)$, one can notice $h$ must lie in $\delta$. Further, by a simple calculation, one can obtain $h=ls$ for any $l\neq0$ and $\alpha=1\otimes s+ls\otimes1+k$.

Summing up two cases above, the solution of (3.2) is $\alpha=0$ or $\alpha=1\otimes s+ls\otimes1+k$ for arbitrary $l, k\in\mathbf{k}$.
$\hfill \blacksquare$
\\

\textbf{Lemma 3.3.} Let $s\in \delta$, $t\in \mathbf{k}$ and $s\neq0$. Then the solution $\alpha \in H\otimes H$ of the following equation
\begin{equation}
  (\alpha\otimes1)(\Delta\otimes \mathrm{id})\alpha=(1\otimes1\otimes s+t)(\mathrm{id}\otimes\Delta)\alpha
\end{equation}
is zero.

\textbf{Proof.} Let $\{\partial_{1}, ..., \partial_{N}\}$ be a basis of $\delta$ and we also consider the corresponding
\emph{PBW} basis given by $\partial^{I}=\partial_{1}^{i_{1}}...\partial_{N}^{i_{n}}/i_{1}!...i_{N}!$. We set $\alpha=\sum_{I}\alpha_{I}\otimes\partial^{I}$, then (3.7) becomes
\begin{equation*}
  \sum_{I}\alpha\Delta(\alpha_{I})\otimes\partial^{I}=\sum_{J+K=I}\alpha_{I}\otimes\partial^{J}\otimes s\partial^{K}+t\sum_{J+K=I}\alpha_{I}\otimes\partial^{J}\otimes\partial^{K}.
\end{equation*}
Let $d$ be the maximal value of $|I|$ for $I$ such that $\alpha_{I}\neq0$. Since $s\neq0$, suppose $\alpha\neq0$, there exist terms in the right $\sum_{|I|=d}\alpha_{I}\otimes1\otimes s\partial^{I}$ lie in $H\otimes \mathbf{k}\otimes F^{d+1}(H)$, which can not be cancelled by other terms whose third tensor factors only contribute lower degree. Then there must have $\alpha=0$.
$\hfill \blacksquare$
\\

\textbf{Lemma 3.4.} Let $s\in \delta$, $t\in \mathbf{k}$ and $s\neq0$. Then the solution $\alpha \in H\otimes H$ of the following equation
\begin{multline}
  (1\otimes s\otimes1+t)(\Delta\otimes \mathrm{id})\alpha-(1\otimes\alpha)(\mathrm{id}\otimes\Delta)(1\otimes s+t)\\
  =(\sigma\otimes \mathrm{id})[(1\otimes s\otimes1+t)(\Delta\otimes \mathrm{id})\alpha-(1\otimes\alpha)(\mathrm{id}\otimes\Delta)(1\otimes s+t)]
\end{multline}
is $\alpha=1\otimes h$ for arbitrary $h\in H$.

\textbf{Proof.} Let $\{\partial_{1}, ..., \partial_{N}\}$ be a basis of $\delta$ and we also consider the corresponding
\emph{PBW} basis given by $\partial^{I}=\partial_{1}^{i_{1}}...\partial_{N}^{i_{n}}/i_{1}!...i_{N}!$. We set $\alpha=\sum_{I}\alpha_{I}\otimes\partial^{I}$, then (3.8) becomes
\begin{align*}
  \sum_{I}[(1\otimes s-s\otimes1)\Delta(\alpha_{I})+(\alpha_{I}\otimes1)(s\otimes1+t)-&(1\otimes\alpha_{I})(1\otimes s+t)]\otimes\partial^{I}\\=&\sum_{I}(1\otimes\alpha_{I}-\alpha_{I}\otimes1)\otimes\partial^{I}s.
\end{align*}
Comparing the degree of third tensor factor of both sides, one can notice $\alpha_{I}$ must lie in $\mathbf{k}$ when $|I|$ reach the maximum value (suppose $|I|_{max}=d$). Therefore, we can rewrite $\alpha=\sum_{|J|<d}\alpha_{J}\otimes\partial^{J}+1\otimes g_{1}$ for some $g_{1}\in F^{d}(H)/F^{d-1}(H)$. Easily verify that $1\otimes g_{1}$ satisfy (3.8) for arbitrary $g_{1}\in H$, then $\sum_{|J|<d}\alpha_{J}\otimes\partial^{J}$ also satisfy (3.8) naturally. By a similar discussion, we have
 \begin{equation*}
\sum_{|J|<d}\alpha_{J}\otimes\partial^{J}=\sum_{|J|<d-1}\alpha_{J}\otimes\partial^{J}+1\otimes g_{2}
\end{equation*}
for some $g_{2}\in F^{d-1}(H)/F^{d-2}(H)$. Through finite times, one can finally obtain $\alpha=1\otimes h$ where $h\in H$ is arbitrary.
$\hfill \blacksquare$
\\

\textbf{Lemma 3.5.} Let $s\in \delta$ and $t, l, k\in \mathbf{k}$. Then the solution $\alpha \in H\otimes H$ of the following equation
\begin{multline}
 (\alpha\otimes1)(\Delta\otimes \mathrm{id})\alpha-(1\otimes1\otimes s+t)(\mathrm{id}\otimes\Delta)\alpha
 =(\sigma\otimes \mathrm{id})[(1\otimes s\otimes1+ls\otimes1\otimes1+k)(\Delta\otimes \mathrm{id})\alpha\\-(1\otimes\alpha)(\mathrm{id}\otimes\Delta)(1\otimes s+ls\otimes1+k)]
\end{multline}
is $\alpha=0$ or $\alpha=t$ or $\alpha=1\otimes s+t$.

\textbf{Proof.} Let $\{\partial_{1}, ..., \partial_{N}\}$ be a basis of $\delta$ and we also consider the corresponding
\emph{PBW} basis given by $\partial^{I}=\partial_{1}^{i_{1}}...\partial_{N}^{i_{n}}/i_{1}!...i_{N}!$. We set $\alpha=\sum_{I}\partial^{I}\otimes\alpha_{I}$, then (3.9) becomes
\begin{multline*}
  \sum_{I,J,K}\partial^{I}\partial^{J}\otimes\alpha_{I}\partial^{K}\otimes\alpha_{J+K}-\sum_{I}\partial^{I}\otimes(1\otimes s+t)\Delta(\alpha_{I})
=\sum_{J+K=I}[s\partial^{J}\otimes\partial^{K}\otimes\alpha_{I}+\partial^{J}\otimes(ls+k)\partial^{K}\otimes\alpha_{I}]\\-(\sum_{I}\partial^{I}\otimes1\otimes\alpha_{I})(s\otimes1\otimes1+1\otimes1\otimes s+1\otimes ls\otimes1+k).
\end{multline*}
Let $d$ be the maximal value of $|I|$ for $I$ such that $\alpha_{I}\neq0$. Then the first summation in the left contains terms whose first tensor factor have degree $2d$ while other terms in both sides have degree of first tensor factor at most $d$. Then there must have $d=0$ and we can write $\alpha=1\otimes g$ for some $g\in H$.
Substituting this into (3.9), we obtain
\begin{equation}
  g\otimes g=(1\otimes s+t)\Delta(g)-1\otimes gs.
\end{equation}
Obviously $g=0$ must be a solution. We suppose $g\neq0$ and $deg(g)=p$. If $p>1$, term $g\otimes g$ has degree $2p$ while other terms contribute degree $p+1$ at most, which gives a contradiction. Therefore, $g$ must lie in $F^{1}(H)$. We can set $g=g'+c$ for some $g'\in\delta$ and $c\in\mathbf{k}$. Substituting this into (3.10), one can easily obtain
\begin{equation*}
  \begin{cases}
c = t \\
g'=0
\end{cases}
\quad or \quad
\begin{cases}
c = t \\
g'=s
\end{cases}.
\end{equation*}

Therefore, we obtain that the solution of (3.9) is $\alpha=0$ or $\alpha=t$ or $\alpha=1\otimes s+t$.
 $\hfill \blacksquare$
\\

Now we consider the semidirect product $\mathcal{P}=He_{1}\oplus He_{2}$, which is generated by two pre-Lie $H$-pseudoalgebras of rank one:
\begin{equation*}
  e_{1}\ast e_{1}=(1\otimes s_{1}+t_{1})\otimes_{H} e_{1}, \quad e_{2}\ast e_{2}=(1\otimes s_{2}+t_{2})\otimes_{H} e_{2},
\end{equation*}
where $s_{1}, s_{2}\in\delta, t_{1}, t_{2}\in\mathbf{k}$.

Suppose
\begin{align*}
  &e_{1}\ast e_{2}=\alpha_{1}\otimes_{H} e_{1}+\alpha_{2}\otimes_{H} e_{2}\\
&e_{2}\ast e_{1}=\beta_{1}\otimes_{H} e_{1}+\beta_{2}\otimes_{H} e_{2}
\end{align*}
for some $\alpha_{1}, \alpha_{2}, \beta_{1}, \beta_{2}\in H\otimes H$. Then by taking $(e_{1},e_{1},e_{2})$, $(e_{1},e_{2},e_{1})$, $(e_{2},e_{1},e_{1})$, $(e_{1},e_{2},e_{2})$, $(e_{2},e_{1},e_{2})$, $(e_{2},e_{2},e_{1})$ into (1.6) and comparing the coefficients of $e_{1}, e_{2}$, we obtain $\mathcal{P}$ is a left pre-Lie $H$-pseudoalgebra if and only if $\alpha_{i}, \beta_{i}$ satisfy the following conditions:
\begin{multline}
  (1\otimes s_{1}\otimes1+t_{1})(\Delta\otimes \mathrm{id})\alpha_{1}-(1\otimes\alpha_{1})(\mathrm{id}\otimes\Delta)(1\otimes s_{1}+t_{1})-(1\otimes\alpha_{2})(\mathrm{id}\otimes\Delta)\alpha_{1}\\=(\sigma\otimes \mathrm{id})[(1\otimes s_{1}\otimes1+t_{1})(\Delta\otimes \mathrm{id})\alpha_{1}-(1\otimes\alpha_{1})(\mathrm{id}\otimes\Delta)(1\otimes s_{1}+t_{1})-(1\otimes\alpha_{2})(\mathrm{id}\otimes\Delta)\alpha_{1}],
\end{multline}
\begin{multline}
  (1\otimes s_{1}\otimes1+t_{1})(\Delta\otimes \mathrm{id})\alpha_{2}-(1\otimes\alpha_{2})(\mathrm{id}\otimes\Delta)\alpha_{2} \\
  =(\sigma\otimes \mathrm{id})[(1\otimes s_{1}\otimes1+t_{1})(\Delta\otimes \mathrm{id})\alpha_{2}-(1\otimes\alpha_{2})(\mathrm{id}\otimes\Delta)\alpha_{2}],
\end{multline}
\begin{multline}
  (\alpha_{1}\otimes1)(\Delta\otimes \mathrm{id})(1\otimes s_{1}+t_{1})+(\alpha_{2}\otimes1)(\Delta\otimes \mathrm{id})\beta_{1}-(1\otimes\beta_{1})(\mathrm{id}\otimes\Delta)(1\otimes s_{1}+t_{1})-(1\otimes\beta_{2})(\mathrm{id}\otimes\Delta)\alpha_{1}
  \\=(\sigma\otimes \mathrm{id})[(\beta_{1}\otimes1)(\Delta\otimes \mathrm{id})(1\otimes s_{1}+t_{1})+(\beta_{2}\otimes1)(\Delta\otimes \mathrm{id})\beta_{1}-(1\otimes1\otimes s_{1}+t_{1})(\mathrm{id}\otimes\Delta)\beta_{1}],
\end{multline}
\begin{multline}
  (\alpha_{2}\otimes1)(\Delta\otimes \mathrm{id})\beta_{2}-(1\otimes\beta_{2})(\mathrm{id}\otimes\Delta)\alpha_{2}=(\sigma\otimes \mathrm{id})[(\beta_{2}\otimes1)(\Delta\otimes \mathrm{id})\beta_{2}-(1\otimes1\otimes s_{1}+t_{1})(\mathrm{id}\otimes\Delta)\beta_{2}],
\end{multline}
\begin{multline}
 (\alpha_{1}\otimes1)(\Delta\otimes \mathrm{id})\alpha_{1}-(1\otimes1\otimes s_{2}+t_{2})(\mathrm{id}\otimes\Delta)\alpha_{1}=(\sigma\otimes \mathrm{id})[(\beta_{1}\otimes1)(\Delta\otimes \mathrm{id})\alpha_{1}-(1\otimes\alpha_{1})(\mathrm{id}\otimes\Delta)\beta_{1}],
\end{multline}
\begin{multline}
  (\alpha_{1}\otimes1)(\Delta\otimes \mathrm{id})\alpha_{2}+(\alpha_{2}\otimes1)(\Delta\otimes \mathrm{id})(1\otimes s_{2}+t_{2})-(1\otimes1\otimes s_{2}+t_{2})(\mathrm{id}\otimes\Delta)\alpha_{2}\\=(\sigma\otimes \mathrm{id})[(\beta_{1}\otimes1)(\Delta\otimes \mathrm{id})\alpha_{2}+(\beta_{2}\otimes1)(\Delta\otimes \mathrm{id})(1\otimes s_{2}+t_{2})-(1\otimes\alpha_{1})(\mathrm{id}\otimes\Delta)\beta_{2}\\-(1\otimes\alpha_{2})(\mathrm{id}\otimes\Delta)(1\otimes s_{2}+t_{2})],
\end{multline}
\begin{multline}
  (1\otimes s_{2}\otimes1+t_{2})(\Delta\otimes \mathrm{id})\beta_{1}-(1\otimes\beta_{1})(\mathrm{id}\otimes\Delta)\beta_{1}=(\sigma\otimes \mathrm{id})[(1\otimes s_{2}\otimes1+t_{2})(\Delta\otimes \mathrm{id})\beta_{1}\\-(1\otimes\beta_{1})(\mathrm{id}\otimes\Delta)\beta_{1}],
\end{multline}
\begin{multline}
 (1\otimes s_{2}\otimes1+t_{2})(\Delta\otimes \mathrm{id})\beta_{2}-(1\otimes\beta_{1})(\mathrm{id}\otimes\Delta)\beta_{2}-(1\otimes\beta_{2})(\mathrm{id}\otimes\Delta)(1\otimes s_{2}+t_{2})\\=(\sigma\otimes \mathrm{id})[(1\otimes s_{2}\otimes1+t_{2})(\Delta\otimes \mathrm{id})\beta_{2}-(1\otimes\beta_{1})(\mathrm{id}\otimes\Delta)\beta_{2}-(1\otimes\beta_{2})(\mathrm{id}\otimes\Delta)(1\otimes s_{2}+t_{2})].
\end{multline}
\\

Combine with results above, we consider the left pre-Lie $H$-pseudoalgebra structures on $\mathcal{P}$ in different situations which are determined by $s_{i}$ and $t_{i}$:

\begin{center}
\begin{tabular}{|c|c|c|c|}
    \hline
    \textbf{Cases} & \textbf{Results} & \textbf{Cases} & \textbf{Results}\\
    \hline
    $s_{i}=t_{i}=0$ & Theorem3.8 & $s_{i}\neq0, t_{i}=0$ & Theorem 3.13\\
    \hline
    $s_{1}\neq0, s_{2}=t_{1}=t_{2}=0$ & Theorem 3.9 & $s_{1}, t_{2}\neq0$, $s_{2}, t_{1}=0$ & Theorem 3.15\\
\hline
    $s_{2}\neq0, s_{1}=t_{1}=t_{2}=0$ & Theorem 3.9 & $s_{2}, t_{1}\neq0$, $s_{1}, t_{2}=0$ & Theorem 3.15\\
\hline
    $t_{1}\neq0, s_{1}=s_{2}=t_{2}=0$ & Theorem 3.10 & $t_{1}=0, s_{1}, s_{2}, t_{2}\neq0$ & Theorem 3.17\\
\hline
    $t_{2}\neq0, s_{1}=s_{2}=t_{1}=0$ & Theorem 3.10 & $t_{2}=0, s_{1}, s_{2}, t_{1}\neq0$ & Theorem 3.17 \\
\hline
$s_{1}, t_{1}\neq0$, $s_{2}, t_{2}=0$ & Theorem 3.9 & $s_{1}=0, s_{2}, t_{1}, t_{2}\neq0$ & Theorem 3.19\\
\hline
    $s_{2}, t_{2}\neq0$, $s_{1}, t_{1}=0$ & Theorem 3.9 & $s_{2}=0, s_{1}, t_{1}, t_{2}\neq0$ & Theorem 3.19\\
\hline
    $t_{i}\neq0, s_{i}=0$ & Theorem 3.11 & $t_{1}, t_{2}, s_{1}, s_{2}\neq0$ & Theorem 3.6\\
\hline
\end{tabular}
\end{center}

\textbf{Theorem 3.6.} If $s_{1}, s_{2}, t_{1}, t_{2}$ are nonzero, that is
\begin{equation*}
  e_{1}\ast e_{1}=(1\otimes s_{1}+t_{1})\otimes_{H} e_{1}, \quad e_{2}\ast e_{2}=(1\otimes s_{2}+t_{2})\otimes_{H} e_{2},
\end{equation*}
 then $\mathcal{P}$ as a left pre-Lie $H$-pseudoalgebra has the following types:

(1) $e_{1}\ast e_{2}=e_{2}\ast e_{1}=0$;

(2) $e_{1}\ast e_{2}=(1\otimes s_{1}+t_{1})\otimes_{H} e_{2}$, $e_{2}\ast e_{1}=(1\otimes s_{1}+t_{1})\otimes_{H} e_{2}$ with $t_{1}s_{2}=t_{2}s_{1}\neq0$;

(3) $e_{1}\ast e_{2}=1\otimes s_{1}\otimes_{H} e_{2}$, $e_{2}\ast e_{1}=1\otimes s_{2}\otimes_{H} e_{1}$ with $[s_{1}, s_{2}]=0$;

(4) $e_{1}\ast e_{2}=(1\otimes s_{1}+t_{1})\otimes_{H} e_{2}$, $e_{2}\ast e_{1}=(1\otimes s_{2}+t_{2})\otimes_{H} e_{1}$ with $t_{1}s_{2}=t_{2}s_{1}\neq0$;

(5) $e_{1}\ast e_{2}=(1\otimes s_{1}+t_{1})\otimes_{H} e_{2}$, $e_{2}\ast e_{1}=1\otimes s_{2}\otimes_{H} e_{1}+t_{1}\otimes1\otimes_{H}e_{2}$ with $[s_{1}, s_{2}]=0$;

(6) $e_{1}\ast e_{2}=(1\otimes s_{1}+2t_{1})\otimes_{H} e_{2}$, $e_{2}\ast e_{1}=(1\otimes s_{2}+\frac{1}{2}t_{2})\otimes_{H} e_{1}+t_{1}\otimes1\otimes_{H}e_{2}$ with $2t_{1}s_{2}=t_{2}s_{1}\neq0$;

(7) $e_{1}\ast e_{2}=t_{2}\otimes1\otimes_{H}e_{1}+1\otimes s_{1}\otimes_{H}e_{2}$, $e_{2}\ast e_{1}=1\otimes s_{2}\otimes_{H}e_{1}+t_{1}\otimes1\otimes_{H}e_{2}$ with $t_{1}s_{2}=t_{2}s_{1}\neq0$;

(8) $e_{1}\ast e_{2}=t_{2}\otimes1\otimes_{H}e_{1}+(1\otimes s_{1}+2t_{1})\otimes_{H} e_{2}$, $e_{2}\ast e_{1}=(1\otimes s_{2}+2t_{2})\otimes_{H} e_{1}+t_{1}\otimes1\otimes_{H}e_{2}$ with $t_{1}s_{2}=t_{2}s_{1}\neq0$;

(9) $e_{1}\ast e_{2}=t_{2}\otimes1\otimes_{H}e_{1}+1\otimes s_{1}\otimes_{H}e_{2}$, $e_{2}\ast e_{1}=(1\otimes s_{2}+t_{2})\otimes_{H} e_{1}$ with $[s_{1}, s_{2}]=0$;

(10) $e_{1}\ast e_{2}=t_{2}\otimes1\otimes_{H}e_{1}+(1\otimes s_{1}+\frac{1}{2}t_{1})\otimes_{H} e_{2}$, $e_{2}\ast e_{1}=(1\otimes s_{2}+2t_{2})\otimes_{H} e_{1}$ with $\frac{1}{2}t_{1}s_{2}=t_{2}s_{1}\neq0$;

(11) $e_{1}\ast e_{2}=(1\otimes s_{2}+t_{2})\otimes_{H} e_{1}$, $e_{2}\ast e_{1}=(1\otimes s_{2}+t_{2})\otimes_{H} e_{1}$ with $t_{1}s_{2}=t_{2}s_{1}\neq0$.

\textbf{Proof.} By equations (3.12), (3.17) and Lemma 3.2, we know
\begin{align*}
  &\alpha_{2}=0 \quad or\quad \alpha_{2}=1\otimes s_{1}+l_{2}s_{1}\otimes1+k_{2},\\
&\beta_{1}=0 \quad or\quad \beta_{1}=1\otimes s_{2}+l_{1}s_{2}\otimes1+k_{1},
\end{align*}
for some $l_{1}, l_{2}, k_{1}, k_{2}\in\mathbf{k}$. We consider in two cases according to $\beta_{1}$:

\textbf{Case I.} When $\beta_{1}=0$.
Then equation (3.15) becomes
\begin{equation*}
  (\alpha_{1}\otimes1)(\Delta\otimes \mathrm{id})\alpha_{1}=(1\otimes1\otimes s_{2}+t_{2})(\mathrm{id}\otimes\Delta)\alpha_{1}.
\end{equation*}
By Lemma 3.3, we know $\alpha_{1}=0$. Further, by equation (3.11), $\alpha_{2}$ is arbitrary. Meanwhile, equation (3.18) becomes
\begin{multline*}
  (1\otimes s_{2}\otimes1+t_{2})(\Delta\otimes \mathrm{id})\beta_{2}-(1\otimes\beta_{2})(\mathrm{id}\otimes\Delta)(1\otimes s_{2}+t_{2})\\=(\sigma\otimes \mathrm{id})[(1\otimes s_{2}\otimes1+t_{2})(\Delta\otimes \mathrm{id})\beta_{2}-(1\otimes\beta_{2})(\mathrm{id}\otimes\Delta)(1\otimes s_{2}+t_{2})].
\end{multline*}
By Lemma 3.4, we obtain $\beta_{2}=1\otimes g$ for some $g\in H$. Then we consider equation (3.14):
\begin{equation}
  (\alpha_{2}\otimes1)(1\otimes g\otimes g)-(1\otimes g\otimes g)(\mathrm{id}\otimes\Delta)\alpha_{2}=g\otimes1\otimes g-(1\otimes1\otimes s_{1}+t_{1})(g_{(1)}\otimes1\otimes g_{(2)}).
\end{equation}
According to $\alpha_{2}$, there are two situations:

If $\alpha_{2}=0$, then we obtain $g\otimes g=(1\otimes s_{1}+t_{1})\Delta(g)$ by (3.19). Suppose $g\neq0$ and $deg(g)=p\geq0$, then there exists term $1\otimes s_{1}g$ lies in $\mathbf{k}\otimes F^{p+1}(H)$ which can not be cancelled by other terms. Thus, there must have $g=0$.

If $\alpha_{2}=1\otimes s_{1}+l_{2}s_{1}\otimes1+k_{2}$, then by equation (3.19), we obtain
 \begin{equation}
g\otimes g+1\otimes gs_{1}=(1\otimes s_{1}+t_{1})\Delta(g).
\end{equation}
 Obviously, $g=0$ is a solution. Suppose $g$ is nonzero and $deg(g)=p\geq0$. If $p>1$, then $2p>p+1$, which implies $g\otimes g\in F^{2p}(H\otimes H)$ can not be cancelled by any other terms. Then $deg(g)\leq1$ and we can write $g=g^{'}+c$ for some $g^{'}\in\delta, c\in\mathbf{k}$. Taking this back to (3.20), we obtain $g=t_{1}$ or $g=s_{1}+t_{1}$.

Summarize discussions above, we have:
\begin{equation*}
\alpha_1 = \beta_1 = 0: \quad
\begin{cases}
1.~\alpha_2=\beta_2 = 0,\\
\alpha_2 = 1 \otimes s_1 + l_2s_1\otimes1 + k_2 :
\begin{cases}
2.~\beta_2 = 0, \\
3.~\beta_2=t_1, \\
4.~\beta_2 = 1\otimes s_1 + t_1.
\end{cases}\\
\end{cases}
\end{equation*}
Easily check that (3.13) and (3.16) will hold when $\alpha_i=\beta_i = 0$. When $\alpha_2 = 1 \otimes s_1 + l_2s_1\otimes1 + k_2$ and $\beta_2 = 1\otimes s_1 + t_1$, one can obtain by (3.16) that $l_{2}=0, k_{2}=t_{1}$ and $s_{i}, t_{i}$ must satisfy $t_{2}s_{1}=t_{1}s_{2}$. When $\alpha_2 = 1 \otimes s_1 + l_2s_1\otimes1 + k_2$, no matter whether $\beta_2=0$ or $\beta_2=t_1$, equation (3.16) will not hold anyway.

Therefore, Case I can contribute to type (1) and type (2).

\textbf{Case II.} When $\beta_{1}=1\otimes s_{2}+l_{1}s_{2}\otimes1+k_{1}$ for some $l_{1}, k_{1}\in\mathbf{k}$.
We set $\alpha_{1}=\sum_{I}\partial^{I}\otimes\alpha_{I}$ and equation (3.15) becomes:
\begin{align*}
  &(\sum_{I} \partial^I \otimes \alpha_I \otimes 1)(\sum_{J+K=I} \partial^J \otimes \partial^K \otimes \alpha_I)- (1\otimes 1\otimes s_2 + t_2)[\sum_{I} \partial^I \otimes \Delta(\alpha_I)]\\
=&(s_2\otimes 1\otimes 1 + 1\otimes l_1 s_2\otimes 1 + k_1)(\sum_{J+K=I} \partial^J \otimes \partial^K \otimes \alpha_I)\\
&-(\sum_{I} \partial^I \otimes 1\otimes \alpha_I)(s_2\otimes 1\otimes 1 + 1\otimes 1\otimes s_2 + 1\otimes l_1 s_2\otimes 1 + k_1).
\end{align*}
Let $d$ be the maximal value of $|I|$ such that $\alpha_{I}\neq0$, then there exists terms in the left whose first tensor factor have degree $2d$ while other terms have degree $d$ at most. Thus, we have $d=0$ and $\alpha_{1}=1\otimes g$ for some $g\in H$. Taking this back to (3.15), we obtain
\begin{equation*}
g\otimes g+1\otimes gs_{2}=(1\otimes s_{2}+t_{2})\Delta(g).
\end{equation*}
 Similar to discussion in Case I, we have $\alpha_{1}=0$ or $\alpha_{1}=t_{2}$ or $\alpha_{1}=1\otimes s_{2}+t_{2}$.

 Then we consider equation (3.11). If $\alpha_{1}=0$, obviously $\alpha_{2}$ is arbitrary. If $\alpha_{1}=t_{2}$, then equation (3.11) becomes $1\otimes\alpha_{2}=(\sigma\otimes\mathrm{id})(1\otimes\alpha_{2})$, one can obtain $\alpha_{2}=0$ or $\alpha_{2}=1\otimes s_{2}+k_{2}$. If $\alpha_{1}=1\otimes s_{2}+t_{2}$, one can verify the only case such that (3.11) holds is $\alpha_{2}=0$.

Finally we discuss $\beta_{2}$ by equation (3.14). If $\alpha_{2}=0$, then equation (3.14) becomes
\begin{equation*}
  (\beta_{2}\otimes1)(\Delta\otimes \mathrm{id})\beta_{2}=(1\otimes1\otimes s_{1}+t_{1})(\mathrm{id}\otimes\Delta)\beta_{2}.
\end{equation*}
By Lemma 3.3, we get $\beta_{2}=0$. If $\alpha_{2}=1\otimes s_{1}+l_{2}s_{1}\otimes1+k_{2}$, then by Lemma 3.5, we have $\beta_{2}=0$ or $\beta_{2}=t_{1}$ or $\beta_{2}=1\otimes s_{1}+t_{1}$. However, one can verify equation (3.18) will not hold anyway when $\beta_{2}=1\otimes s_{1}+t_{1}$.  When $\beta_{2}=t_{1}$, equation (3.18) holds if and only if $\beta_{1}=1\otimes s_{2}+k_{1}$. When $\beta_{2}=0$, equation (3.18) will always hold.

Summarize discussions in Case II, we obtain the following situations:
\begin{align*}
& 1 \begin{cases}
\alpha_1 = 0 \\
\alpha_2 = 0 \\
\beta_1 = 1\otimes s_2+l_1 s_2\otimes 1 + k_1\\
\beta_2 = 0
\end{cases},
\quad
2 \begin{cases}
\alpha_1 = 0 \\
\alpha_2 = 1\otimes s_1+l_2 s_1\otimes 1 + k_2\\
\beta_1 = 1\otimes s_2+l_1 s_2\otimes 1 + k_1\\
\beta_2 = 0
\end{cases},\\
\quad
& 3\begin{cases}
\alpha_1 = 0 \\
\alpha_2 = 1\otimes s_1+l_2 s_1\otimes 1 + k_2\\
\beta_1 = 1\otimes s_2 + k_1\\
\beta_2 = t_1
\end{cases},
\quad
4 \begin{cases}
\alpha_1 = t_2 \\
\alpha_2 = 0 \\
\beta_1 = 1\otimes s_2+l_1 s_2\otimes 1 + k_1\\
\beta_2 = 0
\end{cases},\\
&5 \begin{cases}
\alpha_1 = t_2 \\
\alpha_2 = 1\otimes s_1 + k_2\\
\beta_1 = 1\otimes s_2 + k_1\\
\beta_2 = t_1
\end{cases},
\quad
6 \begin{cases}
\alpha_1 = t_2 \\
\alpha_2 = 1\otimes s_1 + k_2\\
\beta_1 = 1\otimes s_2 \\
\quad \quad +l_1 s_2\otimes 1 + k_1\\
\beta_2 = 0
\end{cases},
\quad
7 \begin{cases}
\alpha_1 = 1\otimes s_2 + t_2\\
\alpha_2 = 0 \\
\beta_1 = 1\otimes s_2 \\
\quad \quad +l_1 s_2\otimes 1 + k_1\\
\beta_2 = 0
\end{cases}.
\end{align*}

 For case 1, one can verify equation (3.13) will not hold anyway.

 For case 2, by taking corresponding $\alpha_{i}, \beta_{i}$ into (3.13) and (3.16), one can obtain
 \begin{equation*}
  \begin{cases}
l_{i}=k_{i}=0\\
[s_{1}, s_{2}]=0
\end{cases}
\quad or \quad
\begin{cases}
l_{i}=0 \\
 k_{1}=t_{2}, k_{2}=t_{1}\\
t_{1}s_{2}=t_{2}s_{1}
\end{cases},
 \end{equation*}
which contributes to type (3) and type (4) respectively.

 For case 3, by taking corresponding $\alpha_{i}, \beta_{i}$ into (3.13) and (3.16), one can obtain
\begin{equation*}
  \begin{cases}
l_{2}=k_{1}=0\\
k_{2}=t_{1}\\
[s_{1}, s_{2}]=0
\end{cases}
\quad or \quad
\begin{cases}
l_{2}=0 \\
 k_{1}=\frac{1}{2}t_{2}, k_{2}=2t_{1}\\
2t_{1}s_{2}=t_{2}s_{1}
\end{cases},
 \end{equation*}
which contributes to type (5) and type (6) respectively.

 For case 4, one can verify equation (3.13) will not hold anyway;

For case 5, by taking corresponding $\alpha_{i}, \beta_{i}$ into (3.13) and (3.16), one can obtain
\begin{equation*}
  \begin{cases}
k_{1}=k_{2}=0\\
t_{1}s_{2}=t_{2}s_{1}
\end{cases}
\quad or \quad
\begin{cases}
k_{1}=2t_{2}, k_{2}=2t_{1} \\
t_{1}s_{2}=t_{2}s_{1}
\end{cases},
 \end{equation*}
which contributes to type (7) and type (8) respectively.

For case 6, by taking corresponding $\alpha_{i}, \beta_{i}$ into (3.13) and (3.16), one can obtain
\begin{equation*}
  \begin{cases}
l_{1}=k_{2}=0\\
k_{1}=t_{2}\\
[s_{1}, s_{2}]=0
\end{cases}
\quad or \quad
\begin{cases}
l_{1}=0 \\
k_{1}=2t_{2}, k_{2}=\frac{1}{2}t_{1}\\
t_{1}s_{2}=2t_{2}s_{1}
\end{cases},
 \end{equation*}
which contributes to type (9) and type (10) respectively.

For case 7, by taking corresponding $\alpha_{i}, \beta_{i}$ into (3.13) and (3.16), one can obtain
$l_{1}=0$, $k_{1}=t_{2}$ and $t_{1}s_{2}=t_{2}s_{1}$,
which contributes to type (11).
$\hfill \blacksquare$
\\

\textbf{Corollary 3.7.} Up to isomorphism, the pre-Lie $H$-pseudoalgebras obtained in Theorem 3.6 can be reduced to the following types:

(i) $e_{1}\ast e_{1}=(1\otimes s_{1}+1)\otimes_{H}e_{1}, e_{1}\ast e_{2}=e_{2}\ast e_{1}=0, e_{2}\ast e_{2}=(1\otimes s_{2}+1)\otimes_{H}e_{2}$, where $s_{1}s_{2}\neq0$;

(ii) $e_{1}\ast e_{1}=(1\otimes s+1)\otimes_{H}e_{1}, e_{1}\ast e_{2}=e_{2}\ast e_{1}=0, e_{2}\ast e_{2}=(1\otimes s+1)\otimes_{H}e_{2}$ for $s\neq0$, which is a special case of type (1);

(iii) $e_{1}\ast e_{1}=(1\otimes s_{1}+1)\otimes_{H}e_{1}, e_{1}\ast e_{2}=1\otimes s_{1}\otimes_{H}e_{2}, e_{2}\ast e_{1}=1\otimes s_{2}\otimes_{H}e_{1}, e_{2}\ast e_{2}=(1\otimes s_{2}+1)\otimes_{H}e_{2}$, where $s_{1}s_{2}\neq0$ and $[s_{1},s_{2}]=0$;

(iv) $e_{1}\ast e_{1}=(1\otimes s+1)\otimes_{H}e_{1}, e_{1}\ast e_{2}=(1\otimes s+1)\otimes_{H}e_{2}, e_{2}\ast e_{1}=e_{2}\ast e_{2}=0$, where $s\neq0$;

(v) $e_{1}\ast e_{1}=(1\otimes s_{1}+1)\otimes_{H}e_{1}, e_{1}\ast e_{2}=(1\otimes s_{1}+1)\otimes_{H}e_{2}, e_{2}\ast e_{1}=1\otimes s_{2}\otimes_{H}e_{1}+1\otimes1\otimes_{H}e_{2}, e_{2}\ast e_{2}=(1\otimes s_{2}+1)\otimes_{H}e_{2}$, where $s_{1}s_{2}\neq0$ and $[s_{1},s_{2}]=0$;

(vi) $e_{1}\ast e_{1}=-1\otimes1\otimes_{H}e_{2}, e_{1}\ast e_{2}=0, e_{2}\ast e_{1}=(1\otimes s+1)\otimes_{H}e_{1}, e_{2}\ast e_{2}=(1\otimes s+2)\otimes_{H}e_{2}$, where $s\neq0$;

(vii) $e_{1}\ast e_{1}=(1\otimes s+1)\otimes_{H}e_{1}, e_{1}\ast e_{2}=1\otimes s\otimes_{H}e_{2}, e_{2}\ast e_{1}=1\otimes1\otimes_{H}e_{2}, e_{2}\ast e_{2}=0$, where $s\neq0$.

\textbf{Proof.}
For type (1) in Theorem 3.6, let $e_{1}'=\frac{e_{1}}{t_{1}}, e_{2}'=\frac{e_{2}}{t_{2}}$ and $s_{1}'=\frac{s_{1}}{t_{1}}, s_{2}'=\frac{s_{2}}{t_{2}}$, then $\{e_{1}', e_{2}'\}$ is an $H$-basis of $\mathcal{P}$ and we have
\begin{equation*}
e_{1}'\ast e_{1}'=(1\otimes s_{1}'+1)\otimes_{H}e_{1}',\quad e_{1}'\ast e_{2}'=e_{2}'\ast e_{1}'=0,\quad e_{2}'\ast e_{2}'=(1\otimes s_{2}'+1)\otimes_{H}e_{2}',
\end{equation*}
where $s_{1}', s_{2}'$ are arbitrary.

For Type (2) in Theorem 3.6, since $t_{1}s_{2}=t_{2}s_{1}\neq0$, let $s\overset{\triangle}{=}\frac{s_{1}}{t_{1}}=\frac{s_{2}}{t_{2}}\neq0$ and $e_{1}'=\frac{e_{1}}{t_{1}}-\frac{e_{2}}{t_{2}}, e_{2}'=\frac{e_{2}}{t_{2}}$, then $\{e_{1}', e_{2}'\}$ is an $H$-basis of $\mathcal{P}$ and we have
\begin{equation*}
e_{1}'\ast e_{1}'=(1\otimes s+1)\otimes_{H}e_{1}', \quad e_{1}'\ast e_{2}'=e_{2}'\ast e_{1}'=0, \quad e_{2}'\ast e_{2}'=(1\otimes s+1)\otimes_{H}e_{2}'.
\end{equation*}

For type (3) in Theorem 3.6, let $e_{1}'=\frac{e_{1}}{t_{1}}, e_{2}'=\frac{e_{2}}{t_{2}}$ and $s_{1}'=\frac{s_{1}}{t_{1}}, s_{2}'=\frac{s_{2}}{t_{2}}$, since $[s_{1}, s_{2}]\neq0$, we have $[s_{1}', s_{2}']\neq0$, then $\{e_{1}', e_{2}'\}$ is an $H$-basis of $\mathcal{P}$ and we have
\begin{align*}
  &e_{1}'\ast e_{1}'=(1\otimes s_{1}+1)\otimes_{H}e_{1}',\quad e_{1}'\ast e_{2}'=1\otimes s_{1}\otimes_{H}e_{2}',\\
 &e_{2}'\ast e_{1}'=1\otimes s_{2}\otimes_{H}e_{1}',\quad e_{2}'\ast e_{2}'=(1\otimes s_{2}+1)\otimes_{H}e_{2}'.
\end{align*}

For type (4) in Theorem 3.6, since $t_{1}s_{2}=t_{2}s_{1}\neq0$, let $s\overset{\triangle}{=}\frac{s_{1}}{t_{1}}=\frac{s_{2}}{t_{2}}\neq0$ and $e_{1}'=\frac{e_{1}}{t_{1}}, e_{2}'=\frac{e_{1}}{t_{1}}-\frac{e_{2}}{t_{2}}$, then $\{e_{1}', e_{2}'\}$ is an $H$-basis of $\mathcal{P}$ and we have
\begin{equation*}
  e_{1}'\ast e_{1}'=(1\otimes s+1)\otimes_{H}e_{1}',\quad e_{1}'\ast e_{2}'=(1\otimes s+1)\otimes_{H}e_{2}',\quad e_{2}'\ast e_{1}'=e_{2}'\ast e_{2}'=0.
\end{equation*}

For type (5) in Theorem 3.6, let $e_{1}'=\frac{e_{1}}{t_{1}}, e_{2}'=\frac{e_{2}}{t_{2}}$ and $s_{1}'=\frac{s_{1}}{t_{1}}, s_{2}'=\frac{s_{2}}{t_{2}}$, since $[s_{1}, s_{2}]\neq0$, we have $[s_{1}', s_{2}']\neq0$, then $\{e_{1}', e_{2}'\}$ is an $H$-basis of $\mathcal{P}$ and we have
\begin{align*}
  &e_{1}'\ast e_{1}'=(1\otimes s_{1}+1)\otimes_{H}e_{1}',\quad e_{1}'\ast e_{2}'=(1\otimes s_{1}+1)\otimes_{H}e_{2}', \\
&e_{2}'\ast e_{1}'=1\otimes s_{2}\otimes_{H}e_{1}'+1\otimes1\otimes_{H}e_{2}',\quad e_{2}'\ast e_{2}'=(1\otimes s_{2}+1)\otimes_{H}e_{2}'.
\end{align*}

For type (6) in Theorem 3.6, since $2t_{1}s_{2}=t_{2}s_{1}\neq0$, let $s\overset{\triangle}{=}\frac{s_{1}}{t_{1}}=\frac{2s_{2}}{t_{2}}\neq0$ and $e_{1}'=\frac{e_{1}}{t_{1}}-2\frac{e_{2}}{t_{2}}, e_{2}'=2\frac{e_{2}}{t_{2}}$, then $\{e_{1}', e_{2}'\}$ is an $H$-basis of $\mathcal{P}$ and we have
\begin{align*}
  &e_{1}'\ast e_{1}'=-1\otimes1\otimes_{H}e_{2}',\quad e_{1}'\ast e_{2}'=0,\\
 &e_{2}'\ast e_{1}'=(1\otimes s+1)\otimes_{H}e_{1}',\quad e_{2}'\ast e_{2}'=(1\otimes s+2)\otimes_{H}e_{2}'.
\end{align*}

For type (7) in Theorem 3.6, since $t_{1}s_{2}=t_{2}s_{1}\neq0$, let $s\overset{\triangle}{=}\frac{s_{1}}{t_{1}}=\frac{s_{2}}{t_{2}}\neq0$ and $e_{1}'=\frac{1}{2}(\frac{e_{1}}{t_{1}}+\frac{e_{2}}{t_{2}}), e_{2}'=\frac{1}{2}(\frac{e_{1}}{t_{1}}-\frac{e_{2}}{t_{2}})$, then $\{e_{1}', e_{2}'\}$ is an $H$-basis of $\mathcal{P}$ and we have
\begin{equation*}
  e_{1}'\ast e_{1}'=(1\otimes s+1)\otimes_{H}e_{1}',\quad e_{1}'\ast e_{2}'=1\otimes s\otimes_{H}e_{2}',\quad
 e_{2}'\ast e_{1}'=1\otimes1\otimes_{H}e_{2}',\quad e_{2}'\ast e_{2}'=0.
\end{equation*}

For type (8) in Theorem 3.6, since $t_{1}s_{2}=t_{2}s_{1}\neq0$, let $s\overset{\triangle}{=}\frac{s_{1}}{t_{1}}=\frac{s_{2}}{t_{2}}\neq0$ and $e_{1}'=\frac{1}{2}(\frac{e_{1}}{t_{1}}+\frac{e_{2}}{t_{2}}), e_{2}'=\frac{1}{2}(\frac{e_{1}}{t_{1}}-\frac{e_{2}}{t_{2}})$, then $\{e_{1}', e_{2}'\}$ is an $H$-basis of $\mathcal{P}$ and we have
\begin{equation*}
  e_{1}'\ast e_{1}'=(1\otimes s+2)\otimes_{H}e_{1}',\quad e_{1}'\ast e_{2}'=(1\otimes s+1)\otimes_{H}e_{2}',\quad e_{2}'\ast e_{1}'=0,\quad e_{2}'\ast e_{2}'=-1\otimes1\otimes_{H}e_{1}',
\end{equation*}
which is actually same as type (6) if we exchange $e_{1}'$ and $e_{2}'$.

For type (9) in Theorem 3.6, let $e_{1}'=\frac{e_{1}}{t_{1}}, e_{2}'=\frac{e_{2}}{t_{2}}$ and $s_{1}'=\frac{s_{1}}{t_{1}}, s_{2}'=\frac{s_{2}}{t_{2}}$, since $[s_{1}, s_{2}]\neq0$, we have $[s_{1}', s_{2}']\neq0$, then $\{e_{1}', e_{2}'\}$ is an $H$-basis of $\mathcal{P}$ and we have
\begin{align*}
  &e_{1}'\ast e_{1}'=(1\otimes s_{1}+1)\otimes_{H}e_{1}',\quad e_{1}'\ast e_{2}'=1\otimes1\otimes_{H}e_{1}'+1\otimes s_{1}\otimes_{H}e_{2}', \\
&e_{2}'\ast e_{1}'=(1\otimes s_{2}+1)\otimes_{H}e_{1}',\quad e_{2}'\ast e_{2}'=(1\otimes s_{2}+1)\otimes_{H}e_{2}',
\end{align*}
which is actually same as type (5) if we exchange $e_{1}'$ and $e_{2}'$.

For type (10) in Theorem 3.6, since $t_{1}s_{2}=2t_{2}s_{1}\neq0$, let $s\overset{\triangle}{=}2\frac{s_{1}}{t_{1}}=\frac{s_{2}}{t_{2}}\neq0$ and $e_{1}'=2\frac{e_{1}}{t_{1}}, e_{2}'=2\frac{e_{1}}{t_{1}}-\frac{e_{2}}{t_{2}}$, then $\{e_{1}', e_{2}'\}$ is an $H$-basis of $\mathcal{P}$ and we have
\begin{equation*}
  e_{1}'\ast e_{1}'=(1\otimes s+2)\otimes_{H}e_{1}',\quad e_{1}'\ast e_{2}'=(1\otimes s+1)\otimes_{H}e_{2}',\quad e_{2}'\ast e_{1}'=0,\quad e_{2}'\ast e_{2}'=-1\otimes1\otimes_{H}e_{1}',
\end{equation*}
which is actually same type (6) if we exchange $e_{1}'$ and $e_{2}'$.

For type (11) in Theorem 3.6, since $t_{1}s_{2}=t_{2}s_{1}\neq0$, let $s\overset{\triangle}{=}\frac{s_{1}}{t_{1}}=\frac{s_{2}}{t_{2}}\neq0$ and $e_{1}'=\frac{e_{1}}{t_{1}}, e_{2}'=\frac{e_{2}}{t_{2}}-\frac{e_{1}}{t_{1}}$, then $\{e_{1}', e_{2}'\}$ is an $H$-basis of $\mathcal{P}$ and we have
\begin{equation*}
e_{1}'\ast e_{1}'=(1\otimes s+1)\otimes_{H}e_{1}', \quad e_{1}'\ast e_{2}'=e_{2}'\ast e_{1}'=0, \quad e_{2}'\ast e_{2}'=(1\otimes s+1)\otimes_{H}e_{2}',
\end{equation*}
which is same as type (2).
$\hfill \blacksquare$
\\

\textbf{Theorem 3.8.} If $s_{1}, s_{2}, t_{1}, t_{2}$ are all zero, that is
\begin{equation*}
  e_{1}\ast e_{1}=0, \quad e_{2}\ast e_{2}=0,
\end{equation*}
then $\mathcal{P}$ as a left pre-Lie $H$-pseudoalgebra has the following types:

(1) $e_{1}\ast e_{2}=0$, $e_{2}\ast e_{1}=0$;

(2) $e_{1}\ast e_{2}=0$, $e_{2}\ast e_{1}=g\otimes1\otimes_{H}e_{1}$ for any nonzero $g\in H$;

(3) $e_{1}\ast e_{2}=h\otimes1\otimes_{H}e_{2}$, $e_{2}\ast e_{1}=0$ for any nonzero $h\in H$.

\textbf{Proof.} By equations (3.12), (3.17) and Lemma 3.1, we know $\beta_{1}=g\otimes1$, $\alpha_{2}=h\otimes1$ for some $g, h\in H$.
Substituting $\alpha_{2}$ in (3.14) and setting $\beta_{2}=\sum_{I}\partial^{I}\otimes\beta_{I}$, we have
\begin{equation*}
  \sum_{J+K=I}(\partial^{K}\otimes h\partial^{J}-\partial^{I}\otimes h)\otimes\beta_{I}=\sum_{I,J,K}\partial^{I}\partial^{J}\otimes\beta_{I}\partial^{K}\otimes\beta_{J+K}.
\end{equation*}
Comparing the degree of first tensor factor in both sides, one can easily obtain $\beta_{2}=0$.

Taking $\beta_{1}=g\otimes1$ into (3.15) and setting $\alpha_{1}=\sum_{I}\partial^{I}\otimes\alpha_{I}$, by a similar discussion, we can also obtain $\alpha_{1}=0$. Then we only need to calculate $g, h$ by (3.13) and (3.16):
\begin{equation*}
  (h\otimes1)\Delta(g)\otimes1=0, \quad (g\otimes1)\Delta(h)\otimes1=0.
\end{equation*}
Obviously, we have $\alpha_{2}=0$ or $\beta_{1}=0$. Thus, there are only three situations:
\begin{equation*}
  1 \begin{cases}
\alpha_1 = 0 \\
\alpha_2 = 0 \\
\beta_1 = 0\\
\beta_2 = 0
\end{cases},
\quad
2 \begin{cases}
\alpha_1 = 0 \\
\alpha_2 = 0\\
\beta_1 = g\otimes1\\
\beta_2 = 0
\end{cases},\\
\quad
3 \begin{cases}
\alpha_1 = 0 \\
\alpha_2 = h\otimes1\\
\beta_1 = 0\\
\beta_2 = 0
\end{cases}
\end{equation*}
for arbitrary nonzero $h, g\in H$. One can easily verify all these three cases satisfy equations (3.11)-(3.18) and then contribute to type (1)-(3) respectively.
$\hfill \blacksquare$
\\

\textbf{Theorem 3.9.} If $s_{1}=t_{1}=0$ and $s_{2}\neq0$, that is
\begin{equation*}
  e_{1}\ast e_{1}=0, \quad e_{2}\ast e_{2}=(1\otimes s+t)\otimes_{H} e_{2}
\end{equation*}
for some $t\in\mathbf{k}$ and nonzero $s\in\delta$, then $\mathcal{P}$ as a left pre-Lie $H$-pseudoalgebra has the following types:

(1) $e_{1}\ast e_{2}=0$, $e_{2}\ast e_{1}=0$;

(2) $e_{1}\ast e_{2}=\alpha\otimes_{H}e_{1}$, $e_{2}\ast e_{1}=(1\otimes s+ls\otimes1+k)\otimes_{H}e_{2}$, where $\alpha\in\{0, t, 1\otimes s+t\}$, $k, l\in\mathbf{k}$ are arbitrary.

\textbf{Proof.} By equation (3.12) and Lemma 3.1, we know $\alpha_{2}=h\otimes1$ for some $h\in H$; by equation (3.17) and Lemma 3.2, we know $\beta_{1}=0$ or $\beta_{1}=1\otimes s+ls\otimes1+k$ for some $l, k\in\mathbf{k}$.

Similar to discussion in Theorem 3.8, we obtain $\beta_{2}=0$ by equation (3.14). Then we discuss in two cases according to $\beta_{1}$:

\textbf{Case I.} When $\beta_{1}=0$. Equation (3.15) becomes
\begin{equation}
   (\alpha_{1}\otimes1)(\Delta\otimes \mathrm{id})\alpha_{1}=(1\otimes1\otimes s+t)(\mathrm{id}\otimes\Delta)\alpha_{1}.
\end{equation}
By Lemma 3.3, we know $\alpha_{1}=0$ and then equation (3.16) becomes
\begin{multline}
  (h\otimes1\otimes1)(1\otimes1\otimes s+t)-(1\otimes1\otimes s+t)(h\otimes1\otimes1)\\=-(\sigma\otimes \mathrm{id})[(1\otimes h\otimes1)(\mathrm{id}\otimes\Delta)(1\otimes s+t)].
\end{multline}
Notice that (3.22) holds if and only if $h=0$. Thus, we have
\begin{equation*}
  \alpha_{1}=\alpha_{2}=\beta_{1}=\beta_{2}=0,
\end{equation*}
which contributes to type (1).

\textbf{Case II.} When $\beta_{1}=1\otimes s+ls\otimes1+k$. By equation (3.15) and Lemma 3.5, we know $\alpha_{1}\in\{0, t, 1\otimes s+t\}$ and we obtain the following situations:
\begin{equation*}
  1 \begin{cases}
\alpha_1 = 0 \\
\alpha_2 = h\otimes1 \\
\beta_1 = 1\otimes s\\
\quad \quad +ls\otimes1+k\\
\beta_2 = 0
\end{cases},
\quad
2 \begin{cases}
\alpha_1 = t \\
\alpha_2 = h\otimes1\\
\beta_1 = 1\otimes s\\
\quad \quad +ls\otimes1+k\\
\beta_2 = 0
\end{cases},\\
\quad
3 \begin{cases}
\alpha_1 = 1\otimes s+t \\
\alpha_2 = h\otimes1\\
\beta_1 = 1\otimes s\\
\quad \quad +ls\otimes1+k\\
\beta_2 = 0
\end{cases}.
\end{equation*}
We only need to calculate the unknown quantities by equations (3.11), (3.13) and (3.16). It is not hard to verify situation 1 (respectively, 2, 3) holds if and only if $h=0$. Therefore, we have
\begin{equation*}
  \begin{cases}
e_{1}\ast e_{2}=\alpha\otimes e_{1}, \quad \alpha\in\{0, t, 1\otimes s+t\}, \\
e_{2}\ast e_{1}= 1\otimes s+ls\otimes1+k\otimes_{H}e_{2}, \quad l, k\in\mathbf{k}~\mbox {are~arbitrary}.
\end{cases}
\end{equation*}
$\hfill \blacksquare$
\\

\textbf{Remark.} The case when $s_{2}=t_{2}=0$ and $s_{1}\neq0$ has no difference with Theorem 3.9 if we exchange $e_{1}$ and $e_{2}$.
\\

\textbf{Theorem 3.10.} If $s_{1}=s_{2}=t_{1}=0$ and $t_{2}\neq0$, that is
\begin{equation*}
  e_{1}\ast e_{1}=0, \quad e_{2}\ast e_{2}=t\otimes1 \otimes_{H}e_{2}
\end{equation*}
for some nonzero $t\in\mathbf{k}$, then $\mathcal{P}$ as a left pre-Lie $H$-pseudoalgebra has the following type:
\begin{equation*}
  e_{1}\ast e_{2}=\alpha\otimes_{H}e_{1}, \quad e_{2}\ast e_{1}=g\otimes1\otimes_{H}e_{1},
\end{equation*}
where $\alpha=0$ or $\alpha=t$, $g\in H$ is arbitrary.

\textbf{Proof.} By equations (3.12), (3.17) and Lemma 3.1, we know $\alpha_{2}=h\otimes1, \beta_{1}=g\otimes1$ for some $h, g\in H$. Similar to discussion in Theorem 3.8, we obtain $\beta_{2}=0$ by equation (3.14). Further, by equations (3.11), (3.13), (3.15) and (3.16), we have
\begin{align}
&(h\otimes1)\Delta(g)\otimes1=0,\\
&\alpha_{1}\Delta(h)+th\otimes1=(1\otimes g)\Delta(h),\\
  &(1\otimes h\otimes1)(\mathrm{id}\otimes\Delta)\alpha_{1}=(\sigma\otimes \mathrm{id})[(1\otimes h\otimes1)(\mathrm{id}\otimes\Delta)\alpha_{1}],\\
&(\alpha_{1}\otimes1)(\Delta\otimes\mathrm{id})\alpha_{1}-t(\mathrm{id}\otimes\Delta)\alpha_{1}=(\sigma\otimes \mathrm{id})[(g\otimes1\otimes1)(\Delta\otimes\mathrm{id})\alpha_{1}-(g\otimes1\otimes1)(1\otimes\alpha_{1})].
\end{align}
Equation (3.23) gives $hg=0$, which induces the following three situations:

\textbf{Case I.} If $g=h=0$, then we only need to consider $(\alpha_{1}\otimes1)(\Delta\otimes\mathrm{id})\alpha_{1}=t(\mathrm{id}\otimes\Delta)\alpha_{1}$. Setting $\alpha_{1}=\sum_{I}\partial^{I}\otimes\alpha_{I}$ and comparing the degree of the first tensor factor in both sides, one can easily notice $\alpha_{1}$ must be $1\otimes\xi$ for some $\xi\in H$. Through a direct calculation, one can obtain $\alpha_{1}=0$ or $\alpha_{1}=t$.

\textbf{Case II.} If $h=0$ and $g\neq0$, we also only need to consider equation (3.26). Similar to case I, one can finally obtain $\alpha_{1}=0$ or $\alpha_{1}=t$.

\textbf{Case III.} If $g=0$ and $h\neq0$, then by equation (3.26), we have $\alpha_{1}=0$ or $\alpha_{1}=t$. By equation (3.24), we have $\alpha_{1}=t$. While in this case, we obtain $th=0$, which is a contradiction.

Summarizing discussion above, we have
\begin{equation*}
  \begin{cases}
e_{1}\ast e_{2}=\alpha\otimes_{H} e_{1}, \quad \alpha\in\{0, t\}, \\
e_{2}\ast e_{1}= g\otimes1\otimes_{H}e_{1}, \quad g\in H~is~arbitrary.
\end{cases}
\end{equation*}
$\hfill \blacksquare$
\\

\textbf{Remark.} The case when $s_{1}=s_{2}=t_{2}=0$ and $t_{1}\neq0$ has no difference with Theorem 3.10 if we exchange $e_{1}$ and $e_{2}$.
\\

\textbf{Theorem 3.11.} If $s_{1}=s_{2}=0$ and $t_{1}, t_{2}\neq0$, that is
\begin{equation*}
  e_{1}\ast e_{1}=t_{1}\otimes1 \otimes_{H}e_{1}, \quad e_{2}\ast e_{2}=t_{2}\otimes1 \otimes_{H}e_{2},
\end{equation*}
 then $\mathcal{P}$ as a left pre-Lie $H$-pseudoalgebra has the following types:

(1) $e_{1}\ast e_{2}=0$, $e_{2}\ast e_{1}=0$;

(2) $e_{1}\ast e_{2}=t_{1}\otimes1 \otimes_{H}e_{2}$, $e_{2}\ast e_{1}=t_{2}\otimes1 \otimes_{H}e_{1}$;

(3) $e_{1}\ast e_{2}=t_{1}\otimes1 \otimes_{H}e_{2}$, $e_{2}\ast e_{1}=t_{1}\otimes1 \otimes_{H}e_{2}$;

(4) $e_{1}\ast e_{2}=2t_{1}\otimes1 \otimes_{H}e_{2}$, $e_{2}\ast e_{1}=\frac{1}{2}t_{2}\otimes1 \otimes_{H}e_{1}+t_{1}\otimes1 \otimes_{H}e_{2}$;

(5) $e_{1}\ast e_{2}=t_{2}\otimes1 \otimes_{H}e_{1}$, $e_{2}\ast e_{1}=t_{2}\otimes1 \otimes_{H}e_{1}$;

(6) $e_{1}\ast e_{2}=t_{2}\otimes1 \otimes_{H}e_{1}$, $e_{2}\ast e_{1}=t_{1}\otimes1 \otimes_{H}e_{2}$;

(7) $e_{1}\ast e_{2}=t_{2}\otimes1 \otimes_{H}e_{1}+\frac{1}{2}t_{1}\otimes1 \otimes_{H}e_{2}$, $e_{2}\ast e_{1}=2t_{2}\otimes1 \otimes_{H}e_{1}$;

(8) $e_{1}\ast e_{2}=t_{2}\otimes1 \otimes_{H}e_{1}+2t_{1}\otimes1 \otimes_{H}e_{2}$, $e_{2}\ast e_{1}=2t_{2}\otimes1 \otimes_{H}e_{1}+t_{1}\otimes1 \otimes_{H}e_{2}$.

\textbf{Proof.} Notice that equation (3.12) can be reduced to
\begin{equation*}
  (1\otimes\alpha_{2})(\mathrm{id}\otimes\Delta)\alpha_{2}=(\sigma\otimes\mathrm{id})[(1\otimes\alpha_{2})(\mathrm{id}\otimes\Delta)\alpha_{2}].
\end{equation*}
By Lemma 3.1, we know $\alpha_{2}=h\otimes1$ for some $h\in H$. Similarly, by equation (3.17), we know $\beta_{1}=g\otimes1$ for some $g\in H$. Setting $\beta_{2}=\sum_{I}\partial^{I}\otimes\beta_{I}$, then equation (3.14) becomes
\begin{equation*}
  \sum_{I}h\partial^{J}\otimes\partial^{K}\otimes\beta_{J+K}-\sum_{I}h\otimes\partial^{I}\otimes\beta_{I}=\sum_{I,J,K}\beta_{I}\partial^{K}\otimes\partial^{I}\partial^{J}\otimes\beta_{J+K}-t_{1}\sum_{I}\beta_{I(1)}\otimes\partial^{I}\otimes\beta_{I(2)}.
\end{equation*}
By comparing the degree of second tensor factors, we know $\beta_{2}=1\otimes u$ for some $u\in H$. Taking this back to (3.14), one can get $\beta_{2}=0$ or $\beta_{2}=t_{1}$. Then equations (3.11), (3.13) and (3.15) are as follows:
\begin{align}
&t_{1}\otimes\alpha_{1}+(1\otimes h\otimes1)(\mathrm{id}\otimes\Delta)\alpha_{1}=(\sigma\otimes\mathrm{id})[t_{1}\otimes\alpha_{1}+(1\otimes h\otimes1)(\mathrm{id}\otimes\Delta)\alpha_{1}]\\
&t_{1}\alpha_{1}\otimes1+(h\otimes1)\Delta(g)\otimes1-1\otimes t_{1}g\otimes1-\beta_{2}(\mathrm{id}\otimes\Delta)\alpha_{1}=\beta_{2}\Delta(g)\otimes1\\
  &(\alpha_{1}\otimes1)(\Delta\otimes\mathrm{id})\alpha_{1}-t_{2}(\mathrm{id}\otimes\Delta)\alpha_{1}=(\sigma\otimes \mathrm{id})[(g\otimes1\otimes1)(\Delta\otimes\mathrm{id})\alpha_{1}-(g\otimes1\otimes1)(1\otimes\alpha_{1})].
\end{align}
Setting $\alpha_{1}=\sum_{I}\partial^{I}\otimes\alpha_{I}$ and comparing the degree of the first tensor factor in both sides in (3.29), we obtain $\alpha_{1}=0$ or $\alpha_{1}=t_{2}$. Then we discuss in two cases according to $\alpha_{1}$:

\textbf{Case I.} When $\alpha_{1}=0$. If $\beta_{2}=0$, then by (3.28), it is easy to know $h, g$ must lie in $\mathbf{k}$. Further, combine with equation (3.16), there are two situations as follows:
\begin{equation*}
  1 \begin{cases}
h = 0 \\
g = 0
\end{cases},
\quad
2 \begin{cases}
h = t_{1} \\
g = t_{2}
\end{cases},
\end{equation*}
which contribute to type (1) and type (2) respectively. If $\beta_{2}=t_{1}$, it is also easy to know $h, g$ must lie in $\mathbf{k}$. Further, combine with equation (3.16), there are two situations as follows:
\begin{equation*}
  3 \begin{cases}
h = t_{1} \\
g = 0
\end{cases},
\quad
4 \begin{cases}
h = 2t_{1} \\
g = \frac{1}{2}t_{2}
\end{cases},
\end{equation*}
which contribute to type (3) and type (4) respectively.

\textbf{Case II.} When $\alpha_{1}=t_{2}$. If $\beta_{2}=0$, then by (3.16), (3.27) and (3.28), $h, g$ have the following possibilities:
\begin{equation*}
  5 \begin{cases}
h = 0\\
g = t_{2}
\end{cases},
\quad
6 \begin{cases}
h = \frac{1}{2}t_{1} \\
g = 2t_{2}
\end{cases},
\end{equation*}
which contribute to type (5) and type (7) respectively. If $\beta_{2}=t_{1}$, then by (3.16), (3.27) and (3.28), $h, g$ have the following possibilities:
\begin{equation*}
  7 \begin{cases}
h = 0\\
g = 0
\end{cases},
\quad
8 \begin{cases}
h = 2t_{1} \\
g = 2t_{2}
\end{cases},
\end{equation*}
which contribute to type (6) and type (8) respectively.
$\hfill \blacksquare$
\\

\textbf{Corollary 3.12.} Up to isomorphism, the pre-Lie pseudoalgebras obtained in Theorem 3.11 can be reduced to the following types:

(i) $e_{1}\ast e_{1}=1\otimes1\otimes_{H}e_{1}, e_{1}\ast e_{2}=e_{2}\ast e_{1}=0, e_{2}\ast e_{2}=1\otimes1\otimes_{H}e_{2}$;

(ii) $e_{1}\ast e_{1}=0, e_{1}\ast e_{2}=0, e_{2}\ast e_{1}=1\otimes1\otimes_{H}e_{1}, e_{2}\ast e_{2}=1\otimes1\otimes_{H}e_{2}$;

(iii) $e_{1}\ast e_{1}=1\otimes1\otimes_{H}2e_{1}, e_{1}\ast e_{2}=1\otimes1\otimes_{H}e_{2}, e_{2}\ast e_{1}=0, e_{2}\ast e_{2}=-1\otimes1\otimes_{H}e_{1}$;

(iv) $e_{1}\ast e_{1}=0, e_{1}\ast e_{2}=1\otimes1\otimes_{H}e_{1}, e_{2}\ast e_{1}=0, e_{2}\ast e_{2}=1\otimes1\otimes_{H}e_{2}$.

\textbf{Proof.} For type (1) in Theorem 3.11, let $e_{1}'=\frac{e_{1}}{t_{1}}, e_{2}'=\frac{e_{2}}{t_{2}}$, then this type is actually isomorphic to (i);

For type (2) in Theorem 3.11, let $e_{1}'=\frac{e_{2}}{t_{2}}-\frac{e_{1}}{t_{1}}, e_{2}'=\frac{e_{2}}{t_{2}}$, then this type is actually isomorphic to (ii);

For type (3) in Theorem 3.11, let $e_{1}'=\frac{e_{1}}{t_{1}}-\frac{e_{2}}{t_{2}}, e_{2}'=\frac{e_{2}}{t_{2}}$, then this type is actually isomorphic to (i);

For type (4) in Theorem 3.11, let $e_{1}'=2\frac{e_{2}}{t_{2}}, e_{2}'=\frac{e_{1}}{t_{1}}-2\frac{e_{2}}{t_{2}}$, then this type is actually isomorphic to (iii);

For type (5) in Theorem 3.11, let $e_{1}'=\frac{e_{1}}{t_{1}}, e_{2}'=\frac{e_{2}}{t_{2}}-\frac{e_{1}}{t_{1}}$, then this type is actually isomorphic to (i);

For type (6) in Theorem 3.11, let $e_{1}'=\frac{e_{1}}{t_{1}}-\frac{e_{2}}{t_{2}}, e_{2}'=\frac{e_{2}}{t_{2}}$, then this type is actually isomorphic to (iv);

For type (7) in Theorem 3.11, let $e_{1}'=2\frac{e_{1}}{t_{1}}, e_{2}'=\frac{e_{2}}{t_{2}}-2\frac{e_{1}}{t_{1}}$, then this type is actually isomorphic to (iii);

For type (8) in Theorem 3.11, let let $e_{1}'=\frac{1}{2}(\frac{e_{1}}{t_{1}}+\frac{e_{2}}{t_{2}}), e_{2}'=\frac{1}{2}(\frac{e_{1}}{t_{1}}-\frac{e_{2}}{t_{2}})$, then this type is actually isomorphic to (iii).
$\hfill \blacksquare$
\\

\textbf{Remark.} Types obtained in Corollary 3.12 are actually 2-dimension pre-Lie algebras, which has been classified and summarized in \cite{BM}.\\

\textbf{Theorem 3.13.} If $t_{1}=t_{2}=0$ and $s_{1}, s_{2}\neq0$, that is
\begin{equation*}
  e_{1}\ast e_{1}=1\otimes s_{1} \otimes_{H}e_{1}, \quad e_{2}\ast e_{2}=1\otimes s_{2} \otimes_{H}e_{2},
\end{equation*}
 then $\mathcal{P}$ as a left pre-Lie $H$-pseudoalgebra has the following types:

(1) $e_{1}\ast e_{2}=0$, $e_{2}\ast e_{1}=0$;

(2) $e_{1}\ast e_{2}=1\otimes s_{1} \otimes_{H}e_{2}$, $e_{2}\ast e_{1}=1\otimes s_{1} \otimes_{H}e_{2}$ with $s_{1}=cs_{2}$ for nonzero $c\in\mathbf{k}$;

(3) $e_{1}\ast e_{2}=1\otimes s_{1} \otimes_{H}e_{2}$, $e_{2}\ast e_{1}=1\otimes s_{2} \otimes_{H}e_{1}$ with $[s_{1}, s_{2}]=0$;

(4) $e_{1}\ast e_{2}=1\otimes s_{2} \otimes_{H}e_{1}$, $e_{2}\ast e_{1}=1\otimes s_{2} \otimes_{H}e_{1}$ with $s_{1}=cs_{2}$ for nonzero $c\in\mathbf{k}$.

\textbf{Proof.} By equations (3.12), (3.17) and Lemma 3.2, we know
\begin{align*}
  &\alpha_{2}=0 \quad or\quad \alpha_{2}=1\otimes s_{1}+l_{2}s_{1}\otimes1+k_{2};\\
&\beta_{1}=0 \quad or\quad \beta_{1}=1\otimes s_{2}+l_{1}s_{2}\otimes1+k_{1}.
\end{align*}
We consider in two cases according to $\beta_{1}$:

\textbf{Case I.} When $\beta_{1}=0$. By equation (3.15) and Lemma 3.3, we know $\alpha_{1}=0$; by equation (3.18) and Lemma 3.4, we know $\beta_{2}=1\otimes g$ for some $g\in H$. Then we further consider $g$ in equation (3.14). If $\alpha_{2}=0$, we obtain $g\otimes g=(1\otimes s_{1})\Delta(g)$, which implies $g=0$. If $\alpha_{2}=1\otimes s_{1}+l_{2}s_{1}\otimes1+k_{2}$ for some $k_{2}, l_{2}\in\mathbf{k}$, we obtain $g\otimes g+1\otimes gs_{1}=(1\otimes s_{1})\Delta(g)$. Through a simple discussion like (3.20), one can obtain $g=0$ or $g=s_{1}$. Therefore, we obtain the following situations:
\begin{equation*}
  1 \begin{cases}
\alpha_1 = 0 \\
\alpha_2 = 0 \\
\beta_1 = 0\\
\beta_2 = 0
\end{cases},
\quad
2 \begin{cases}
\alpha_1 = 0 \\
\alpha_2 = 1\otimes s_{1}+l_{2}s_{1}\otimes1+k_{2} \\
\beta_1 = 0\\
\beta_2 = 0
\end{cases},
\quad
 3 \begin{cases}
\alpha_1 = 0 \\
\alpha_2 = 1\otimes s_{1}+l_{2}s_{1}\otimes1+k_{2} \\
\beta_1 = 0\\
\beta_2 = 1\otimes s_{1}
\end{cases}.
\end{equation*}
Case 1 is obviously hold for all equations. Notice that Case 2 does not hold in equation (3.16) anyway. In case 3, one can verify equations (3.11)-(3.18) hold if and only if $l_{2}=k_{2}=0$ and $s_{1}=cs_{2}$ for nonzero $c\in\mathbf{k}$. Thus, situations above contribute to type (1) and type (2) respectively.

\textbf{Case II.} When $\beta_{1}=1\otimes s_{2}+l_{1}s_{2}\otimes1+k_{1}$. By equation (3.15) and Lemma 3.5, we know $\alpha_{1}=0$ or $\alpha_{1}=1\otimes s_{2}$.

When $\alpha_{1}=0$, then by equation (3.11), $\alpha_{2}$ is arbitrary. When $\alpha_{1}=1\otimes s_{2}$, one can verify (3.11) will not hold anyway if $\alpha_{2}=1\otimes s_{1}+l_{2}s_{1}\otimes1+k_{2}$. Then $\alpha_{2}$ can only be zero in this situation.

Further, when $\alpha_{2}=0$, by equation (3.14) and Lemma 3.3, we know $\beta_{2}=0$. When $\alpha_{2}=1\otimes s_{1}+l_{2}s_{1}\otimes1+k_{2}$, then by equation (3.14) and Lemma 3.5, we know $\beta_{2}=0$ or $\beta_{2}=1\otimes s_{1}$. However, $\beta_{2}=1\otimes s_{1}$ must be abandoned since equation (3.18) will not hold in this case anyway.

According to discussion above, we obtain the following situations:
\begin{equation*}
  1 \begin{cases}
\alpha_1 = 0 \\
\alpha_2 = 0 \\
\beta_1 = 1\otimes s_{2}\\
\quad \quad +l_{1}s_{2}\otimes1+k_{1}\\
\beta_2 = 0
\end{cases},
2 \begin{cases}
\alpha_1 = 0 \\
\alpha_2 = 1\otimes s_{1}+l_{2}s_{1}\otimes1+k_{2} \\
\beta_1 = 1\otimes s_{2}+l_{1}s_{2}\otimes1+k_{1}\\
\beta_2 = 0
\end{cases},
 3 \begin{cases}
\alpha_1 = 1\otimes s_{2} \\
\alpha_2 = 0 \\
\beta_1 = 1\otimes s_{2}\\
\quad \quad +l_{1}s_{2}\otimes1+k_{1}\\
\beta_2 = 0
\end{cases}.
\end{equation*}
By verifying equations (3.13) and (3.16), case 1 will not hold anyway; case 2 holds if and only if $l_{i}=k_{i}=0$ and $[s_{1}, s_{2}]=0$; case 3 holds if and only if $l_{1}=k_{1}=0$ and $s_{1}=cs_{2}$ for nonzero $c\in\mathbf{k}$, which contribute to type (3) and type (4) respectively.
$\hfill \blacksquare$
\\

\textbf{Corollary 3.14.} Up to isomorphism, the pre-Lie pseudoalgebras obtained in Theorem 3.13 can be reduced to the following types:

(i) $e_{1}\ast e_{1}=1\otimes s_{1}\otimes_{H}e_{1}, e_{1}\ast e_{2}=e_{2}\ast e_{1}=0, e_{2}\ast e_{2}=1\otimes s_{2}\otimes_{H}e_{2}$, where $s_{1}s_{2}\neq0$;

(ii) $e_{1}\ast e_{1}=1\otimes s\otimes_{H}e_{1}, e_{1}\ast e_{2}=e_{2}\ast e_{1}=0, e_{2}\ast e_{2}=1\otimes s\otimes_{H}e_{2}$, which is a special case of type (1);

(iii) $e_{1}\ast e_{1}=1\otimes s_{1}\otimes_{H}e_{1}, e_{1}\ast e_{2}=1\otimes s_{1}\otimes_{H}e_{2}, e_{2}\ast e_{1}=1\otimes s_{2}\otimes_{H}e_{1}, e_{2}\ast e_{2}=1\otimes s_{2}\otimes_{H}e_{2}$, where $s_{1}s_{2}\neq0$ and $[s_{1},s_{2}]=0$.

\textbf{Proof.} The proof is similar to Corollary 3.7.
$\hfill \blacksquare$
\\

\textbf{Theorem 3.15.} If $s_{1}=t_{2}=0$ and $s_{2}, t_{1}$ are nonzero, that is
\begin{equation*}
  e_{1}\ast e_{1}=t\otimes1 \otimes_{H}e_{1}, \quad e_{2}\ast e_{2}=1\otimes s \otimes_{H}e_{2}
\end{equation*}
for some nonzero $s\in\delta$ and nonzero $t\in\mathbf{k}$, then $\mathcal{P}$ as a left pre-Lie $H$-pseudoalgebra has the following types:

(1) $e_{1}\ast e_{2}=0$, $e_{2}\ast e_{1}=0$;

(2) $e_{1}\ast e_{2}=0$, $e_{2}\ast e_{1}=1\otimes s\otimes_{H}e_{1}$;

(3) $e_{1}\ast e_{2}=t\otimes1\otimes_{H}e_{2}$, $e_{2}\ast e_{1}=1\otimes s\otimes_{H}e_{1}+t\otimes1\otimes_{H}e_{2}$;

(4) $e_{1}\ast e_{2}=1\otimes s\otimes_{H}e_{1}$, $e_{2}\ast e_{1}=(1\otimes s+s\otimes1)\otimes_{H}e_{1}$.

\textbf{Proof.}  By equation (3.12) and Lemma 3.1, we know $\alpha_{2}=h\otimes1$ for some $h\in H$. By equation (3.17) and Lemma 3.2, we know $\beta_{1}=0$ or $\beta_{1}=1\otimes s+ls\otimes1+k$ for some $l, k\in\mathbf{k}$. Setting $\beta_{2}=\sum_{I}\partial^{I}\otimes\beta_{I}$, then equation (3.14) becomes
\begin{equation*}
  \sum_{I}h\partial^{J}\otimes\partial^{K}\otimes\beta_{J+K}-\sum_{I}h\otimes\partial^{I}\otimes\beta_{I}=\sum_{I,J,K}\beta_{I}\partial^{K}\otimes\partial^{I}\partial^{J}\otimes\beta_{J+K}-t_{1}\sum_{I}\beta_{I(1)}\otimes\partial^{I}\otimes\beta_{I(2)}.
\end{equation*}
By comparing the degree of second tensor factors, we know $\beta_{2}=1\otimes u$ for some $u\in H$. Taking this back to (3.14), one can obtain $\beta_{2}=0$ or $\beta_{2}=t$. Then we consider $\alpha_{1}$ by equation (3.15) in two cases:

\textbf{Case I.} When $\beta_{1}=0$. Then we have $(\alpha_{1}\otimes1)(\Delta\otimes \mathrm{id})\alpha_{1}=(1\otimes1\otimes s+t)(\mathrm{id}\otimes\Delta)\alpha_{1}$, which implies $\alpha_{1}=0$ by Lemma 3.3.

\textbf{Case II.} When $\beta_{1}=1\otimes s+ls\otimes1+k$ for some $l, k\in\mathbf{k}$. Immediately by Lemma 3.5, we obtain $\alpha_{1}=0$ or $\alpha_{1}=1\otimes s$.

According to discussion above, we obtain the following situations:
\begin{align*}
   &1 \begin{cases}
\alpha_1 = 0 \\
\alpha_2 = h\otimes1 \\
\beta_1 =0\\
\beta_2 = 0
\end{cases},
\quad
 2 \begin{cases}
\alpha_1 = 0 \\
\alpha_2 = h\otimes1 \\
\beta_1 =0\\
\beta_2 = t
\end{cases},
\quad
 3 \begin{cases}
\alpha_1 = 0 \\
\alpha_2 = h\otimes1 \\
\beta_1 = 1\otimes s+ls\otimes1+k\\
\beta_2 = 0
\end{cases}\\
\quad
&4 \begin{cases}
\alpha_1 = 0 \\
\alpha_2 = h\otimes1 \\
\beta_1 = 1\otimes s+ls\otimes1+k \\
\beta_2 = t
\end{cases},
\quad
5 \begin{cases}
\alpha_1 = 1\otimes s \\
\alpha_2 = h\otimes1 \\
\beta_1 = 1\otimes s+ls\otimes1+k\\
\beta_2 = 0
\end{cases},
\quad
6 \begin{cases}
\alpha_1 = 1\otimes s \\
\alpha_2 = h\otimes1 \\
\beta_1 = 1\otimes s \\
\quad \quad +ls\otimes1+k\\
\beta_2 = t
\end{cases}.
\end{align*}

We verify these situations in equations (3.11), (3.13), and (3.16).

For case 1, by taking  corresponding $\alpha_{i}, \beta_{i}$ into equation (3.16), one can obtain $\alpha_{2}=0$, which contributes to type (1).

For case 2, one can verify equation (3.16) will not hold anyway.

For case 3, by taking  corresponding $\alpha_{i}, \beta_{i}$ into equations (3.13) and (3.16), one can obtain $l=k=h=0$, which contributes to type (2).

For case 4, by taking  corresponding $\alpha_{i}, \beta_{i}$ into equations (3.13) and (3.16), one can obtain $l=k=0$, $h=t$, which contributes to type (3).

For case 5, by taking  corresponding $\alpha_{i}, \beta_{i}$ into equations (3.11) and (3.13), one can obtain $k=0, l=1, h=0$, which contributes to type (4).

For case 6, one can obtain $h=0$ by equation (3.11) while equation (3.13) will not hold anyway in this situation.
$\hfill \blacksquare$
\\

\textbf{Remark.} The case when $s_{2}=t_{1}=0$ and $s_{1}, t_{2}$ are nonzero has no difference with Theorem 3.15 if we exchange $e_{1}$ and $e_{2}$.
\\

\textbf{Corollary 3.16.} Up to isomorphism, the pre-Lie pseudoalgebras obtained in Theorem 3.15 can be reduced to the following types:

(i) $e_{1}\ast e_{1}=1\otimes 1\otimes_{H}e_{1}, e_{1}\ast e_{2}=e_{2}\ast e_{1}=0, e_{2}\ast e_{2}=1\otimes s\otimes_{H}e_{2}$, where $s\neq0$;

(ii) $e_{1}\ast e_{1}=1\otimes 1\otimes_{H}e_{1}, e_{1}\ast e_{2}=0, e_{2}\ast e_{1}=1\otimes s\otimes_{H}e_{1}, e_{2}\ast e_{2}=1\otimes s\otimes_{H}e_{2}$, where $s\neq0$;

(iii) $e_{1}\ast e_{1}=1\otimes 1\otimes_{H}e_{1}, e_{1}\ast e_{2}=1\otimes 1\otimes_{H}e_{2}, e_{2}\ast e_{1}=1\otimes s\otimes_{H}e_{1}+1\otimes 1\otimes_{H}e_{2}, e_{2}\ast e_{2}=1\otimes s\otimes_{H}e_{2}$, where $s\neq0$.

\textbf{Proof.} We only discuss type (4) in Theorem 3.15: let $e_{1}'=\frac{e_{1}}{t}, e_{2}'=e_{2}-\frac{se_{1}}{t}$, then $\{e_{1}', e_{2}'\}$ is also an $H$-basis and one can calculate that
\begin{equation*}
  e_{1}'\ast e_{1}'=1\otimes 1\otimes_{H}e_{1}',\quad e_{1}'\ast e_{2}'=0,\quad e_{2}'\ast e_{1}'=1\otimes s\otimes_{H}e_{1}',\quad e_{2}'\ast e_{2}'=1\otimes s\otimes_{H}e_{2}',
\end{equation*}
which contributes to type (iii).
$\hfill \blacksquare$
\\

\textbf{Theorem 3.17.} If $t_{1}=0$ and $s_{1}, s_{2}$ are nonzero, that is
\begin{equation*}
  e_{1}\ast e_{1}=1\otimes s_{1} \otimes_{H}e_{1}, \quad e_{2}\ast e_{2}=(1\otimes s_{2}+t)\otimes_{H}e_{2}
\end{equation*}
for some nonzero $t\in\mathbf{k}$, then $\mathcal{P}$ as a left pre-Lie $H$-pseudoalgebra has the following types:

(1) $e_{1}\ast e_{2}=0$, $e_{2}\ast e_{1}=0$;

(2) $e_{1}\ast e_{2}=1\otimes s_{1}\otimes_{H}e_{2}$, $e_{2}\ast e_{1}=1\otimes s_{2}\otimes_{H}e_{1}$ with $[s_{1}, s_{2}]=0$;

(3) $e_{1}\ast e_{2}=t\otimes1\otimes_{H}e_{1}+1\otimes s_{1}\otimes_{H}e_{2}$, $e_{2}\ast e_{1}=(1\otimes s_{2}+t)\otimes_{H}e_{1}$ with $[s_{1}, s_{2}]=0$.

\textbf{Proof.} By equations (3.12), (3.17) and Lemma 3.2, we know
\begin{align*}
  &\alpha_{2}=0 \quad or\quad \alpha_{2}=1\otimes s_{1}+l_{2}s_{1}\otimes1+k_{2};\\
&\beta_{1}=0 \quad or\quad \beta_{1}=1\otimes s_{2}+l_{1}s_{2}\otimes1+k_{1}.
\end{align*}
We consider in two cases according to $\beta_{1}$:

\textbf{Case I.} When $\beta_{1}=0$. By equation (3.15) and Lemma 3.3, we know $\alpha_{1}=0$. By equation (3.18) and Lemma 3.4, we know $\beta_{2}=1\otimes g$ for some $g\in H$. If $\alpha_{2}=0$, by equation (3.14), we have $g\otimes g=(1\otimes s_{1})\Delta(g)$, which implies $\beta_{2}=0$. If $\alpha_{2}=1\otimes s_{1}+l_{2}s_{1}\otimes1+k_{2}$, then by equation (3.14) and Lemma 3.5, we know $\beta_{2}=0$ or $\beta_{2}=1\otimes s_{1}$. Thus, we obtain the following situations:
\begin{equation*}
  1 \begin{cases}
\alpha_1 = 0 \\
\alpha_2 = 0 \\
\beta_1 = 0\\
\beta_2 = 0
\end{cases},
\quad
2 \begin{cases}
\alpha_1 = 0 \\
\alpha_2 = 1\otimes s_{1}+l_{2}s_{1}\otimes1+k_{2} \\
\beta_1 = 0\\
\beta_2 = 0
\end{cases},
\quad
 3 \begin{cases}
\alpha_1 = 0 \\
\alpha_2 = 1\otimes s_{1}+l_{2}s_{1}\otimes1+k_{2} \\
\beta_1 = 0\\
\beta_2 = 1\otimes s_{1}
\end{cases}.
\end{equation*}
Case 1 obviously holds for all equations. One can verify both case 2 and case 3 will not hold in equation (3.16) anyway.

\textbf{Case II.} When $\beta_{1}=1\otimes s_{2}+l_{1}s_{2}\otimes1+k_{1}$. This case is not different from Theorem 3.6 if we let $t_{1}=0$. Thus, there are five situations as follows:
\begin{align*}
&1 \begin{cases}
\alpha_1 = 0 \\
\alpha_2 = 0 \\
\beta_1 = 1\otimes s_2\\
\quad \quad +l_1 s_2\otimes 1 + k_1\\
\beta_2 = 0
\end{cases},
\quad
2 \begin{cases}
\alpha_1 = 0 \\
\alpha_2 = 1\otimes s_1+l_2 s_1\otimes 1 + k_2\\
\beta_1 = 1\otimes s_2+l_1 s_2\otimes 1 + k_1\\
\beta_2 = 0
\end{cases},
\quad
3 \begin{cases}
\alpha_1 = t\otimes 1 \\
\alpha_2 = 0 \\
\beta_1 = 1\otimes s_2\\
\quad \quad +l_1 s_2\otimes 1 + k_1\\
\beta_2 = 0
\end{cases},\\
&4 \begin{cases}
\alpha_1 = t\otimes 1 \\
\alpha_2 = 1\otimes s_1 + k_2\otimes 1\\
\beta_1 = 1\otimes s_2+l_1 s_2\otimes 1 + k_1\\
\beta_2 = 0
\end{cases},
\quad
5 \begin{cases}
\alpha_1 = 1\otimes s_2 + t\otimes 1\\
\alpha_2 = 0 \\
\beta_1 = 1\otimes s_2+l_1 s_2\otimes 1 + k_1\\
\beta_2 = 0
\end{cases}.
\end{align*}
For cases 1, 3 and 5, one can verify equations (3.13) and (3.16) do not hold anyway.

For case 2, by equations (3.13) and (3.16), one can obtain $k_{i}=l_{i}=0$ and there must have $[s_{1}, s_{2}]=0$, which contributes to type (2).

For case 4, by equations (3.13) and (3.16), one can obtain $l_{1}=k_{2}=0$, $k_{1}=t$ and there must have $[s_{1}, s_{2}]=0$, which contributes to type (3).
$\hfill \blacksquare$
\\

\textbf{Remark.} The case when $t_{2}=0$ and $t_{1}, s_{1}, s_{2}$ are nonzero has no difference with Theorem 3.17 if we exchange $e_{1}$ and $e_{2}$.
\\

\textbf{Corollary 3.18.} Up to isomorphism, the pre-Lie pseudoalgebras obtained in Theorem 3.17 can be reduced to the following types:

(i) $e_{1}\ast e_{1}=1\otimes s_{1}\otimes_{H}e_{1}, e_{1}\ast e_{2}=e_{2}\ast e_{1}=0, e_{2}\ast e_{2}=(1\otimes s_{2}+1)\otimes_{H}e_{2}$, where $s_{1}, s_{2}$ are arbitrary and $s_{1}s_{2}\neq0$;

(ii) $e_{1}\ast e_{1}=1\otimes s_{1}\otimes_{H}e_{1}, e_{1}\ast e_{2}=1\otimes s_{1}\otimes_{H}e_{2}, e_{2}\ast e_{1}=1\otimes s_{2}\otimes_{H}e_{1}, e_{2}\ast e_{2}=(1\otimes s_{2}+1)\otimes_{H}e_{2}$, where $s_{1}s_{2}\neq0$ and $[s_{1}, s_{2}]=0$;

(iii) $e_{1}\ast e_{1}=1\otimes s_{1}\otimes_{H}e_{1}, e_{1}\ast e_{2}=1\otimes1\otimes_{H}e_{1}+1\otimes s_{1}\otimes_{H}e_{2}, e_{2}\ast e_{1}=(1\otimes s_{2}+1)\otimes_{H}e_{1}, e_{2}\ast e_{2}=(1\otimes s_{2}+1)\otimes_{H}e_{2}$, where $s_{1}s_{2}\neq0$ and $[s_{1}, s_{2}]=0$.

\textbf{Proof.} The proof is similar to Corollary 3.7.
$\hfill \blacksquare$
\\

\textbf{Theorem 3.19.} If $s_{1}=0$ and $s_{2}\neq0$, that is
\begin{equation*}
  e_{1}\ast e_{1}=t_{1}\otimes1 \otimes_{H}e_{1}, \quad e_{2}\ast e_{2}=(1\otimes s+t_{2})\otimes_{H}e_{2},
\end{equation*}
where $t_{1}, t_{2}\in\mathbf{k}$, $s\in\delta$ are nonzero, then $\mathcal{P}$ as a left pre-Lie $H$-pseudoalgebra has the following types:

(1) $e_{1}\ast e_{2}=0$, $e_{2}\ast e_{1}=0$;

(2) $e_{1}\ast e_{2}=0$, $e_{2}\ast e_{1}=1\otimes s\otimes_{H}e_{1}$;

(3) $e_{1}\ast e_{2}=t_{1}\otimes1\otimes_{H}e_{2}$, $e_{2}\ast e_{1}=1\otimes s\otimes_{H}e_{1}+t_{1}\otimes1\otimes_{H}e_{2}$;

(4) $e_{1}\ast e_{2}=t_{2}\otimes1\otimes_{H}e_{1}$, $e_{2}\ast e_{1}=(1\otimes s+t_{2})\otimes_{H}e_{1}$;

(5) $e_{1}\ast e_{2}=(1\otimes s+t_{2})\otimes_{H}e_{1}$, $e_{2}\ast e_{1}=(1\otimes s+s\otimes1+t_{2})\otimes_{H}e_{1}$.

\textbf{Proof.} By equation (3.12) and Lemma 3.1, we know $\alpha_{2}=h\otimes1$ for some $h\in H$; by equation (3.17) and Lemma 3.2, we know $\beta_{1}=0$ or $\beta_{1}=1\otimes s+ls\otimes1+k$ for some $k, l\in\mathbf{k}$. Then we consider in two cases according to $\beta_{1}$:

\textbf{Case I.} When $\beta_{1}=0$. By equation (3.15) and Lemma 3.3, we know $\alpha_{1}=0$; by equation (3.18) and Lemma 3.4 we know $\beta_{2}=1\otimes g$ for some $g\in H$. Then by equation (3.14), we have $g\otimes g=t_{1}\Delta(g)$, which implies $g=0$ or $g=t_{1}$. Thus, we obtain the following situations:
\begin{equation*}
  1 \begin{cases}
\alpha_1 = 0 \\
\alpha_2 = h\otimes1 \\
\beta_1 = 0\\
\beta_2 = 0
\end{cases}
\quad and \quad
2 \begin{cases}
\alpha_1 = 0 \\
\alpha_2 = h\otimes1\\
\beta_1 = 0\\
\beta_2 = t_{1}
\end{cases}.
\end{equation*}
By verifying equations (3.13) and (3.16), one can notice equation (3.16) does not hold anyway in case 2. Case 1 holds for (3.13) and (3.16) if and only if $h=0$, which contributes to type (1).

\textbf{Case II.} When $\beta_{1}=1\otimes s+ls\otimes1+k$. By equation (3.15) and Lemma 3.5, we obtain $\alpha_{1}\in\{0, t_{2}, 1\otimes s+t_{2}\}$. By equation (3.11), if $\alpha_{1}=0$, then $\alpha_{2}$ is arbitrary; if $\alpha_{1}=t_{2}$, then $\alpha_{2}\in\mathbf{k}$; if $\alpha_{1}=1\otimes s+t_{2}$, then $\alpha_{2}$ must be zero.

Next we consider $\beta_{2}$ by equation (3.14). Let $\beta_{2}=\sum_{I}\partial^{I}\otimes\beta_{I}$, then we have
\begin{equation*}
  \sum_{I}h\partial^{J}\otimes\partial^{K}\otimes\beta_{I}-\sum_{I}h\otimes\partial^{I}\otimes\beta_{I}=\sum_{I,J,K}\beta_{I}\partial^{K}\otimes\partial^{I}\partial^{J}\otimes\beta_{J+K}-t_{1}\sum_{I}\beta_{I(1)}\otimes\partial^{I}\otimes\beta_{I(2)}.
\end{equation*}
Let $d$ be the maximal value of $|I|$ for $I$ such that $\alpha_{I}\neq0$. Notice that there exists terms in the right whose second tensor factor have degree $2d$ while other terms have degree at most $d$. Then we have $d=0$ and $\beta_{2}=1\otimes g$ for some $g\in H$. Further, one can obtain $g=0$ or $g=t_{1}$ by taking $1\otimes g$ back to (3.14).

Finally, through equation (3.18), we obtain $\beta_{1}$ is arbitrary when $\beta_{2}=0$ and $\beta_{1}=1\otimes s+k$ for arbitrary $k$ when $\beta_{2}=t_{1}$. Summarizing discussion above, we have the following situations:
\begin{align*}
  &1 \begin{cases}
\alpha_1 = 0 \\
\alpha_2 = h\otimes1 \\
\beta_1 = 1\otimes s+l s\otimes 1 + k\\
\beta_2 = 0
\end{cases},
\quad
2 \begin{cases}
\alpha_1 = 0 \\
\alpha_2 = h\otimes1\\
\beta_1 = 1\otimes s+k\\
\beta_2 = t_{1}
\end{cases},
\quad
3 \begin{cases}
\alpha_1 = t_2 \\
\alpha_2 = h\in\mathbf{k} \\
\beta_1 = 1\otimes s+ls\otimes 1 + k\\
\beta_2 = 0
\end{cases},\\
&4 \begin{cases}
\alpha_1 = t_2\\
\alpha_2 = h\in\mathbf{k}\\
\beta_1 = 1\otimes s+k\\
\beta_2 = t_{1}
\end{cases},
\quad
5 \begin{cases}
\alpha_1 = 1\otimes s + t_2\\
\alpha_2 = 0 \\
\beta_1 = 1\otimes s+ls\otimes 1 + k\\
\beta_2 = 0
\end{cases},
\quad
6 \begin{cases}
\alpha_1 = 1\otimes s + t_2\\
\alpha_2 = 0 \\
\beta_1 = 1\otimes s+k\\
\beta_2 = t_{1}
\end{cases}.
\end{align*}

For case 1, by by taking corresponding $\alpha_{i}, \beta_{i}$ into equations (3.13) and (3.16), one can obtain $h=l=k=0$, which contributes to type (2);

For case 2, by by taking corresponding $\alpha_{i}, \beta_{i}$ into equations (3.13) and (3.16), one can obtain $h=t_{1}, k=0$, which contributes to type (3);

For case 3, by by taking corresponding $\alpha_{i}, \beta_{i}$ into equations (3.13) and (3.16), one can obtain $h=l=0$, $k=t_{2}$, which contributes to type (4);

For case 4, by by taking corresponding $\alpha_{i}, \beta_{i}$ into equation (3.13), one can obtain $k=0$, $h=t_{1}$, while this gives a contradiction with (3.16);

For case 5, by by taking corresponding $\alpha_{i}, \beta_{i}$ into equations (3.13) and (3.16), one can obtain $l=1$, $k=t_{2}$, which contributes to type (5);

For case 6, one can verify equation (3.13) will not hold anyway.
$\hfill \blacksquare$
\\

\textbf{Remark.} The case when $s_{2}=0$ and $s_{1}, t_{1}, t_{2}$ are zero has no difference with Theorem 3.19 if we exchange $e_{1}$ and $e_{2}$.
\\

\textbf{Corollary 3.20.} Up to isomorphism, the pre-Lie pseudoalgebras obtained in Theorem 3.19 can be reduced to the following types:

(i) $e_{1}\ast e_{1}=1\otimes 1\otimes_{H}e_{1}, e_{1}\ast e_{2}=e_{2}\ast e_{1}=0, e_{2}\ast e_{2}=(1\otimes s+1)\otimes_{H}e_{2}$, where $s\neq0$;

(ii) $e_{1}\ast e_{1}=1\otimes 1\otimes_{H}e_{1}, e_{1}\ast e_{2}=0, e_{2}\ast e_{1}=1\otimes s\otimes_{H}e_{1}, e_{2}\ast e_{2}=(1\otimes s+1)\otimes_{H}e_{2}$, where $s\neq0$;

(iii) $e_{1}\ast e_{1}=1\otimes 1\otimes_{H}e_{1}, e_{1}\ast e_{2}=1\otimes 1\otimes_{H}e_{2}, e_{2}\ast e_{1}=1\otimes s\otimes_{H}e_{1}+1\otimes 1\otimes_{H}e_{2}, e_{2}\ast e_{2}=(1\otimes s+1)\otimes_{H}e_{2}$, where $s\neq0$.

\textbf{Proof.} We only discuss type (5) in Theorem 3.19: let $e_{1}'=\frac{e_{1}}{t_{1}}, e_{2}'=\frac{e_{2}}{t_{2}}-\frac{(s+t_{2})e_{1}}{t_{1}t_{2}}$, then $\{e_{1}', e_{2}'\}$ is also an $H$-basis and one can calculate that
\begin{equation*}
  e_{1}'\ast e_{1}'=1\otimes 1\otimes_{H}e_{1}',\quad e_{1}'\ast e_{2}'=0,\quad e_{2}'\ast e_{1}'=1\otimes s\otimes_{H}e_{1}',\quad e_{2}'\ast e_{2}'=(1\otimes s+1)\otimes_{H}e_{2}',
\end{equation*}
which contributes to type (ii).
$\hfill \blacksquare$
\\

\section*{4. Associative $H$-pseudoalgebras}
\def\theequation{4. \arabic{equation}}
	\setcounter{equation} {0} \hskip\parindent

In this section, we consider the associative pseudoalgebra structure on $\mathcal{P}$ and further we will point out which left pre-Lie pseudoalgebras obtained in Section 3 are associative. We still discussion under the assumption that $H$ is an universal enveloping algebra.
\\

Let $He$ be a rank one $H$-pseudoalgebra. If $He$ is associative, then there have $(e\ast e)\ast e=e\ast (e\ast e)$, which implies
\begin{equation}
(\alpha\otimes1)(\Delta\otimes\mathrm{id})\alpha=(1\otimes\alpha)(\mathrm{id}\otimes\Delta)\alpha,
\end{equation}
if we let $e\ast e=\alpha\otimes_{H}e$ for some $\alpha\in H\otimes H$.

We can write $\alpha=\sum_{I}\alpha_{I}\otimes\partial^{I}$ for some $\alpha_{I}\in H$, then (4.1) becomes
\begin{equation*}
  \sum_{I}\alpha\Delta(\alpha_{I})\otimes\partial^{I}=\sum_{I,J,K}\alpha_{J+K}\otimes\alpha_{I}\partial^{J}\otimes\partial^{I}\partial^{K}.
\end{equation*}
One can easily obtain that there must have $|I|=0$ and $\alpha=g\otimes1$ for some $g\in H$. Further, by substituting this into (4.1), we have $(g\otimes1)\Delta(g)=g\otimes g$, which implies $\alpha=0$ or $\alpha\in\mathbf{k}$. Thus, we have the following result:
\\

\textbf{Theorem 4.1.} The $He$ is an rank one associative $H$-pseudoalgebra if and only if $e\ast e=t\otimes1\otimes_{H}e$ for $t\in\mathbf{k}$. Further, through a simple basis transformation, there are two types up to isomorphism:

(1) $e\ast e=0$,

(2) $e\ast e=1\otimes1\otimes_{H}e$.\\

Then we consider the rank two pseudoalgebra $\mathcal{P}=He_{1}\oplus He_{2}$. By discussion above, if $\mathcal{P}$ is associative, there must have $e_{1}\ast e_{1}=t_{1}\otimes1\otimes_{H}e_{1}$, $e_{2}\ast e_{2}=t_{2}\otimes1\otimes_{H}e_{2}$ for $t_{1}, t_{2}\in\mathbf{k}$.

Suppose
\begin{align*}
  &e_{1}\ast e_{2}=\alpha_{1}\otimes_{H} e_{1}+\alpha_{2}\otimes_{H} e_{2}\\
&e_{2}\ast e_{1}=\beta_{1}\otimes_{H} e_{1}+\beta_{2}\otimes_{H} e_{2}
\end{align*}
for some $\alpha_{1}, \alpha_{2}, \beta_{1}, \beta_{2}\in H\otimes H$. Then by $H$-bilinearity and taking $(e_{1},e_{1},e_{2})$, $(e_{1},e_{2},e_{1})$, $(e_{2},e_{1},e_{1})$, $(e_{1},e_{2},e_{2})$, $(e_{2},e_{1},e_{2})$, $(e_{2},e_{2},e_{1})$ into (1.6) to compare the coefficients of $e_{1}, e_{2}$, we obtain $\mathcal{P}$ is an associative $H$-pseudoalgebra if and only if $\alpha_{i}, \beta_{i}$ satisfy the following conditions:
\begin{align}
  &t_{1}(\Delta\otimes\mathrm{id})\alpha_{1}=t_{1}(1\otimes\alpha_{1})+(1\otimes\alpha_{2})(\mathrm{id}\otimes\Delta)\alpha_{1},\\
&t_{1}(\Delta\otimes\mathrm{id})\alpha_{2}=(1\otimes\alpha_{2})(\mathrm{id}\otimes\Delta)\alpha_{2},\\
&t_{1}(\alpha_{1}\otimes1)+(\alpha_{2}\otimes1)(\Delta\otimes\mathrm{id})\beta_{1}=t_{1}(1\otimes\beta_{1})+(1\otimes\beta_{2})(\mathrm{id}\otimes\Delta)\alpha_{1},\\
&(\alpha_{2}\otimes1)(\Delta\otimes\mathrm{id})\beta_{2}=(1\otimes\beta_{2})(\mathrm{id}\otimes\Delta)\alpha_{2},\\
&t_{1}(\beta_{1}\otimes1)+(\beta_{2}\otimes1)(\Delta\otimes\mathrm{id})\beta_{1}=t_{1}(\mathrm{id}\otimes\Delta)\beta_{1},\\
&t_{1}(\mathrm{id}\otimes\Delta)\beta_{2}=(\beta_{2}\otimes1)(\Delta\otimes\mathrm{id})\beta_{2},\\
&(\alpha_{1}\otimes1)(\Delta\otimes\mathrm{id})\alpha_{1}=t_{2}(\mathrm{id}\otimes\Delta)\alpha_{1},\\
&t_{2}(\mathrm{id}\otimes\Delta)\alpha_{2}=(\alpha_{1}\otimes1)(\Delta\otimes\mathrm{id})\alpha_{2}+t_{2}(\alpha_{2}\otimes1),\\
&(\beta_{1}\otimes1)(\Delta\otimes\mathrm{id})\alpha_{1}=(1\otimes\alpha_{1})(\mathrm{id}\otimes\Delta)\beta_{1},\\
&(\beta_{1}\otimes1)(\Delta\otimes\mathrm{id})\alpha_{2}+t_{2}(\beta_{2}\otimes1)=(1\otimes\alpha_{1})(\mathrm{id}\otimes\Delta)\beta_{2}+t_{2}(1\otimes\alpha_{2}),\\
&t_{2}(\Delta\otimes\mathrm{id})\beta_{1}=(1\otimes\beta_{1})(\mathrm{id}\otimes\Delta)\beta_{1},\\
&t_{2}(\Delta\otimes\mathrm{id})\beta_{2}=(1\otimes\beta_{1})(\mathrm{id}\otimes\Delta)\beta_{2}+t_{2}(1\otimes\beta_{2}).
\end{align}
\\

We still consider the structure in different situations which are determined by $t_{1}, t_{2}$:\\

\textbf{Theorem 4.2.} If
 $e_{1}\ast e_{1}=0, e_{2}\ast e_{2}=0$,
 then $\mathcal{P}$ is associative if and only if $e_{1}\ast e_{2}=e_{2}\ast e_{1}=0$.

\textbf{Proof.} From equation (4.3), we know $\alpha_{2}=0$; from equation (4.7), we know $\beta_{2}=0$; from equation (4.8), we know $\alpha_{1}=0$; from equation (4.12), we know $\beta_{1}=0$. Obviously, other equations will hold naturally.
$\hfill \blacksquare$
\\

\textbf{Theorem 4.3.} If $e_{1}\ast e_{1}=0, e_{2}\ast e_{2}=1\otimes1 \otimes_{H}e_{2}$, then $\mathcal{P}$ is associative if and only if

(1) $e_{1}\ast e_{2}=e_{2}\ast e_{1}=0$;

(2) $e_{1}\ast e_{2}=0, e_{2}\ast e_{1}=1\otimes1 \otimes_{H}e_{1}$;

(3) $e_{1}\ast e_{2}=1\otimes1 \otimes_{H}e_{1}, e_{2}\ast e_{1}=0$;

(4) $e_{1}\ast e_{2}=1\otimes1 \otimes_{H}e_{1}, e_{2}\ast e_{1}=1\otimes1 \otimes_{H}e_{1}$.

\textbf{Proof.} From equation (4.3), we know $\alpha_{2}=0$; from equation (4.7), we know $\beta_{2}=0$; from equation (4.8), we know $\alpha_{1}=0$ or $\alpha_{1}=1$; from equation (4.12), we know $\beta_{1}=0$ or $\beta_{1}=1$. Therefore, we obtain four situations as follows:
\begin{align*}
  &1 \begin{cases}
e_{1}\ast e_{2}=0 \\
e_{2}\ast e_{1}=0
\end{cases},
2 \begin{cases}
e_{1}\ast e_{2}=0 \\
e_{2}\ast e_{1}=1\otimes1 \otimes_{H}e_{1}
\end{cases},\\
&3 \begin{cases}
e_{1}\ast e_{2}=1\otimes1 \otimes_{H}e_{1} \\
e_{2}\ast e_{1}=0
\end{cases},
4 \begin{cases}
e_{1}\ast e_{2}=1\otimes1 \otimes_{H}e_{1}\\
e_{2}\ast e_{1}=1\otimes1 \otimes_{H}e_{1}
\end{cases}.
\end{align*}
One can verify all these situations satisfy equations (4.2)-(4.13) and contribute to type (1) to type (4) respectively.
$\hfill \blacksquare$
\\

\textbf{Remark.} The case when $e_{1}\ast e_{1}=1\otimes1 \otimes_{H}e_{1}, e_{2}\ast e_{2}=0$ has no difference with Theorem 4.3 if we exchange $e_{1}$ and $e_{2}$.
\\

\textbf{Theorem 4.4.} If
  $e_{1}\ast e_{1}=t_{1}\otimes1 \otimes_{H}e_{1}, e_{2}\ast e_{2}=t_{2}\otimes1 \otimes_{H}e_{2}$
for some nonzero $t_{1}, t_{2}\in\mathbf{k}$, then $\mathcal{P}$ is associative if and only if

(1) $e_{1}\ast e_{2}=e_{2}\ast e_{1}=0$;

(2) $e_{1}\ast e_{2}=t_{2}\otimes1 \otimes_{H}e_{1}, e_{2}\ast e_{1}=t_{2}\otimes1 \otimes_{H}e_{1}$;

(3) $e_{1}\ast e_{2}=t_{1}\otimes1 \otimes_{H}e_{2}, e_{2}\ast e_{1}=t_{2}\otimes1 \otimes_{H}e_{1}$;

(4) $e_{1}\ast e_{2}=t_{2}\otimes1 \otimes_{H}e_{1}, e_{2}\ast e_{1}=t_{1}\otimes1 \otimes_{H}e_{2}$;

(5) $e_{1}\ast e_{2}=t_{1}\otimes1 \otimes_{H}e_{2}, e_{2}\ast e_{1}=t_{1}\otimes1 \otimes_{H}e_{2}$.

\textbf{Proof.} From equation (4.3), we know $\alpha_{2}=0$ or $\alpha_{2}=t_{1}$; from equation (4.7), we know $\beta_{2}=0$ or $\beta_{2}=t_{1}$; from equation (4.8), we know $\alpha_{1}=0$ or $\alpha_{1}=t_{2}$; from equation (4.12), we know $\beta_{1}=0$ or $\beta_{1}=t_{2}$. Therefore, we obtain the following situations:
\begin{align*}
  &1 \begin{cases}
e_{1}\ast e_{2}=0 \\
e_{2}\ast e_{1}=0
\end{cases},
2 \begin{cases}
e_{1}\ast e_{2}=t_{1}\otimes1 \otimes_{H}e_{2} \\
e_{2}\ast e_{1}=0
\end{cases},
3 \begin{cases}
e_{1}\ast e_{2}=0 \\
e_{2}\ast e_{1}=t_{2}\otimes1 \otimes_{H}e_{1}
\end{cases},\\
&4 \begin{cases}
e_{1}\ast e_{2}=t_{1}\otimes1 \otimes_{H}e_{2}\\
e_{2}\ast e_{1}=t_{2}\otimes1 \otimes_{H}e_{1}
\end{cases},
5 \begin{cases}
e_{1}\ast e_{2}=0\\
e_{2}\ast e_{1}=t_{1}\otimes1 \otimes_{H}e_{2}
\end{cases},
6 \begin{cases}
e_{1}\ast e_{2}=t_{1}\otimes1 \otimes_{H}e_{2}\\
e_{2}\ast e_{1}=t_{1}\otimes1 \otimes_{H}e_{2}
\end{cases},\\
&7\begin{cases}
e_{1}\ast e_{2}=0\\
e_{2}\ast e_{1}=t_{2}\otimes1 \otimes_{H}e_{1}+t_{1}\otimes1 \otimes_{H}e_{2}
\end{cases},
8 \begin{cases}
e_{1}\ast e_{2}=t_{1}\otimes1 \otimes_{H}e_{2}\\
e_{2}\ast e_{1}=t_{2}\otimes1 \otimes_{H}e_{1}+t_{1}\otimes1 \otimes_{H}e_{2}
\end{cases},\\
&9 \begin{cases}
e_{1}\ast e_{2}=t_{2}\otimes1 \otimes_{H}e_{1} \\
e_{2}\ast e_{1}=0
\end{cases},
10 \begin{cases}
e_{1}\ast e_{2}=t_{2}\otimes1 \otimes_{H}e_{1} \\
e_{2}\ast e_{1}=t_{2}\otimes1 \otimes_{H}e_{1}
\end{cases},
11 \begin{cases}
e_{1}\ast e_{2}=t_{2}\otimes1 \otimes_{H}e_{1} \\
e_{2}\ast e_{1}=t_{1}\otimes1 \otimes_{H}e_{2}
\end{cases},\\
&12 \begin{cases}
e_{1}\ast e_{2}=t_{2}\otimes1 \otimes_{H}e_{1}+t_{1}\otimes1 \otimes_{H}e_{2} \\
e_{2}\ast e_{1}=0
\end{cases},
13 \begin{cases}
e_{1}\ast e_{2}=t_{2}\otimes1 \otimes_{H}e_{1}+t_{1}\otimes1 \otimes_{H}e_{2} \\
e_{2}\ast e_{1}=t_{2}\otimes1 \otimes_{H}e_{1}
\end{cases},\\
&14 \begin{cases}
e_{1}\ast e_{2}=t_{2}\otimes1 \otimes_{H}e_{1}+t_{1}\otimes1 \otimes_{H}e_{2} \\
e_{2}\ast e_{1}=t_{1}\otimes1 \otimes_{H}e_{2}
\end{cases},
15 \begin{cases}
e_{1}\ast e_{2}=t_{2}\otimes1 \otimes_{H}e_{1} \\
e_{2}\ast e_{1}=t_{2}\otimes1 \otimes_{H}e_{1}+t_{1}\otimes1 \otimes_{H}e_{2}
\end{cases},\\
&16 \begin{cases}
e_{1}\ast e_{2}=t_{2}\otimes1 \otimes_{H}e_{1}+t_{1}\otimes1 \otimes_{H}e_{2} \\
e_{2}\ast e_{1}=t_{2}\otimes1 \otimes_{H}e_{1}+t_{1}\otimes1 \otimes_{H}e_{2}
\end{cases}.
\end{align*}
Notice that equation (4.2) does not hold in case 12, 13, 14 and 16; equation (4.4) does not hold in case 3, 7 and 9; equation (4.11) does not hold in case 2 and 5; equation (4.13) does not hold in case 8 and 15. The rest of the situations will always satisfy equations (4.2)-(4.13) and contribute to type (1) to type (5) respectively.
$\hfill \blacksquare$
\\

\textbf{Corollary 4.5.} Up to isomorphism, the associative $H$-pseudoalgebras obtained in Theorem 4.4 can be reduced to the following types:

(i) $e_{1}\ast e_{1}=1\otimes1 \otimes_{H}e_{1}, e_{1}\ast e_{2}=e_{2}\ast e_{1}=0, e_{2}\ast e_{2}=1\otimes1 \otimes_{H}e_{2}$;

(ii) $e_{1}\ast e_{1}=0, e_{1}\ast e_{2}=0, e_{2}\ast e_{1}=1\otimes1 \otimes_{H}e_{1}, e_{2}\ast e_{2}=1\otimes1 \otimes_{H}e_{2}$;

(iii) $e_{1}\ast e_{1}=0, e_{1}\ast e_{2}=1\otimes1 \otimes_{H}e_{1}, e_{2}\ast e_{1}=0, e_{2}\ast e_{2}=1\otimes1 \otimes_{H}e_{2}$.

\textbf{Proof.} The proof is similar to Corollary 3.7.\\

Notice that all types obtained in Theorem 4.2 to Theorem 4.4 can be found in Theorem 3.8, 3.10 and 3.11 respectively, which means $\mathcal{P}$ as a left pre-Lie $H$-pseudoalgebra is associative if and only if $\mathcal{P}$ is one of the types of Theorem 4.2 to Theorem 4.4.\\

\section*{DATA AVAILABILITY}

Data sharing is not applicable to this article as no new data were created or analyzed in
this study.

\end{document}